\documentclass{amsart}
\usepackage{graphicx}
\usepackage{url}

\usepackage[all]{xy} 
\input xy 

\xyoption{all}
\xyoption{2cell} 
\xyoption{v2}

\DeclareMathOperator{\tr}{tr}

\DeclareMathOperator{\spec}{spec}
\DeclareMathOperator{\reg}{reg} 
\DeclareMathOperator{\sgn}{sgn}

\theoremstyle{plain}
\newtheorem{theorem}{Theorem}[section]
\newtheorem{lemma}[theorem]{Lemma}
\newtheorem{corollary}[theorem]{Corollary}
\newtheorem{proposition}[theorem]{Proposition}

\theoremstyle{definition}
\newtheorem{definition}[theorem]{Definition}

\theoremstyle{remark}

\numberwithin{equation}{section}
\numberwithin{figure}{section}

\renewcommand{\cL}{{\mathcal L}}

\newcommand{\CC}{{\mathbb C}}
\newcommand{\RR}{{\mathbb R}}
\newcommand{\ZZ}{{\mathbb Z}}

\newcommand{\TT}{{\mathbb T}} 

\renewcommand{\a}{\alpha}

\newcommand{\Res}{\mathrm{Res}}

\begin{document}

\title[Bundle gerbe]{The basic bundle gerbe \\ on unitary groups}

 \author[M. Murray]{Michael Murray}
 \address[Michael Murray]
 {School of Mathematical Sciences\\
 University of Adelaide\\
 Adelaide, SA 5005 \\
 Australia}
 \email{michael.murray@adelaide.edu.au}

\author[D. Stevenson]{Danny Stevenson}
\address[Danny Stevenson]
{Fachbereich Mathematik\\
Bundesstrasse 55 \\
Universit\"{a}t Hamburg\\
D-20146 Hamburg\\
Germany}
\email{stevenson@math.uni-hamburg.de} 

 \thanks{The first author acknowledges the support of the Australian
Research Council and thanks Alan Carey and Mathai Varghese for useful conversations. The second author thanks the University of Adelaide 
for hospitality during which time part of this research was carried out and acknowledges 
support from the Collaborative Research Center 676 `Particles, strings and the early universe'.}


\begin{abstract}
We consider the construction of the basic bundle gerbe on $SU(n)$ introduced by Meinrenken and show that it extends
to a range of groups with unitary actions on a Hilbert space including $U(n)$, $\TT^n$ and $U_p(H)$, the Banach Lie group of unitaries differing
from the identity by an element of a Schatten ideal. In all these cases we give   an explicit connection and curving on the basic bundle gerbe and calculate
the real  Dixmier-Douady class. Extensive use is made of the holomorphic functional calculus for operators on a Hilbert space. 
\end{abstract}
\maketitle
\tableofcontents

\section{Introduction}
Gerbes were introduced by Giraud \cite{Gir} to study  non-abelian cohomology.   Brylinski popularised them in his book \cite{Bry}, in particular
the case of interest here, which is gerbes with band the sheaf of smooth $U(1)$ valued functions.  To every such gerbe 
is associated a characteristic class in $H^3(M,\ZZ)$, the Dixmier-Douady class of the gerbe.  Equivalence classes of 
gerbes on $M$ are, through this characteristic class, in bijective correspondence with 
$H^3(M,\ZZ)$.    Therefore it is natural to look for gerbes on manifolds with 
non-zero degree three integer cohomology and one important  example of such a manifold is a compact, simple  and simply connected Lie group.  

Recall that if $G$ is a compact, simple and simply connected 
Lie group then $H^3(G,\ZZ) = \ZZ$ and there is a canonical closed three-form $\nu$ on $G$ --- 
the \emph{basic three-form} (see for instance \cite{PS}).  $\nu$ is a de Rham representative 
for the image in real cohomology of the generator of $H^3(G,\ZZ)$.  The three-form $\nu$  
is given by 
$$ 
\nu = \frac{1}{12}\langle \theta_L,[\theta_L,\theta_L]\rangle 
$$ 
where $\theta_L$ is the left Maurer-Cartan 1-form on $G$ and 
$\langle \, ,\, \rangle$ is the \emph{basic} inner product on 
$\mathfrak{g}$ \cite{PS}.  
In the case where $G = SU(n)$ the 
basic three-form is 
\begin{equation} 
\label{eq: basic 3-form on SU(n)} 
-\frac{1}{24\pi^2}\tr(g^{-1}dg)^3.  
\end{equation}
Although the unitary group $U(n)$ is not simply connected, we still have the isomorphism 
$H^3(U(n),\ZZ) = \ZZ$.  The three-form~\eqref{eq: basic 3-form on SU(n)} on $SU(n)$ is clearly the 
restriction of a closed 3-form defined on $U(n)$.  This three-form is the image in real cohomology 
of the generator of $H^3(U(n),\ZZ)$ --- we will refer to it as the basic three-form on $U(n)$.  

The basic three-form $\nu$ was exploited to good effect in Witten's paper \cite{Witten} on 
WZW models. Witten considered a non-linear sigma-model 
in which the fields of the theory were smooth maps $g\colon \Sigma\to G$ from 
a compact Riemann surface $\Sigma$ to a compact, simple and simply connected Lie group 
$G$.  In constructing a conformally invariant action for this sigma-model Witten 
was lead to consider the Wess-Zumino term 
$$ 
S_{WZ}(g) = \int_B \tilde{g}^*\nu. 
$$ 
Here $B$ is a 3-manifold with boundary $\Sigma$ and $\tilde{g}\colon B\to G$ is an 
extension of $g$ to $B$.  The question arises as to whether the Wess-Zumino 
term is well defined.  It turns out that under these topological assumptions on $G$, one 
can always find such an extension of $g$, and, due to the integrality property of the basic three-form $\nu$, 
$\exp(2\pi iS_{WZ}(g))$ is well defined.  It is natural to wonder if it is possible 
to remove the topological assumptions on $G$ and make sense of this 
action when $G$ is a non-simply connected group.   The theory of gerbes provides a valuable 
way of thinking about this problem (see for 
example \cite{Schweigert-Waldorf}).  One can interpret the action $\exp(2\pi i S_{WZ}(g))$ 
as the holonomy over $\Sigma$ of a canonically defined gerbe on $G$ --- the \emph{basic gerbe}.  
This basic gerbe on $G$ can be constructed even when the group $G$ is not 
simply connected.   Since 
the holonomy of a gerbe on a manifold can be defined irrespective of whether the 
manifold is simply connected or not, we see that by \emph{defining} the action 
to be the holonomy of the basic gerbe over $\Sigma$, we can remove this 
topological assumption on the group $G$.     

There have been a number of constructions of gerbes and bundle gerbes on a Lie group $G$ 
in the literature since Brylinksi's book \cite{Bry} appeared and we review them briefly to put 
the results of this paper into perspective. Indeed the  first such construction appeared in  \cite{Bry} (and later in \cite{BryMac}); it involved the path-fibration 
of $G$ and thus was  inherently infinite-dimensional. It was pointed out in \cite{Bry} that 
it would be interesting to have a finite dimensional construction.   

The notion of bundle gerbe was introduced by the first author in  \cite{Mur96Bundle-gerbes}.  The relationship of bundle gerbes with gerbes is analogous to that between 
line bundles thought of as fibrations and line bundles thought of as locally free sheaves of modules. Bundles gerbes correspond to fibrations of groupoids whereas 
gerbes involve sheaves of groupoids.  In \cite{Mur96Bundle-gerbes} the 
tautological bundle gerbe was introduced.  This was a bundle gerbe associated to any integral, closed three-form on a 2-connected manifold $M$. This implicitly includes the case of 
compact, simple, simply-connected Lie groups which was discussed more explicitly in \cite{CMW}.  Again these constructions are infinite-dimensional and related to the path fibration.  There is a simple
way of defining this bundle gerbe using the so-called lifting bundle gerbe described in \cite{Mur96Bundle-gerbes}.  
The path-fibration over $G$ is a principal bundle with structure group $\Omega G$, the group of based loops in $G$, and there is a well known central extension $\widehat{\Omega G}$ of 
$\Omega G$ by the circle (see \cite{PS}) which one can use to form a bundle gerbe.  This bundle gerbe measures   the obstruction to lifting the path-fibration 
to a bundle with structure group $\widehat{\Omega G}$.

The next construction, due to Brylinski \cite{Bry1} (see also \cite{Brytalk}), involves the Weyl map
\begin{equation*}
\begin{array}{ccc}
G/T \times T & \to &G \\
(gT, t) & \mapsto & gtg^{-1}
\end{array}
\end{equation*}
A gerbe was defined on $G/T\times T$ using a `cup product' construction involving line bundles on $G/T$ and functions on $T$.   
It was shown, using some delicate sheaf arguments, that this  gerbe pushes forward via the Weyl map to a  gerbe on $G$.  Brylinski notes that this 
construction gives an equivariant gerbe for
the conjugation action of $G$ on $G$.  This construction of Brylinski seems to be 
the most general construction to date, however the geometry of this gerbe has not been explored in full detail.  

Following this construction of Brylinski's was a construction of Gawedzki and Reis \cite{Gawedzki-Reis} for 
the case when $G = SU(n)$.  The case of quotients of $SU(n)$ by finite subgroups of the centre 
was also treated.  The methods used in this construction involved ideas from representation theory.  Gawedzki and Reis 
also defined a connection and curving on their bundle gerbe.  This work was followed shortly afterwards by a paper of 
Meinrenken \cite{Meinrenken}.  This gave a definitive treatment of the case where $G$ was an arbitrary compact, simple and 
simply connected Lie group.  This construction was also representation theoretic in nature and involved the 
structure of the sets of regular and singular elements of $G$.  Meinrenken's paper also gave an extensive 
discussion of equivariant bundle gerbes; the basic bundle gerbe constructed in the paper was shown to be equivariant and 
equipped with an equivariant connection and curving.    A simpler construction 
of a local bundle gerbe in the sense of Chatterjee-Hitchin was also given for the case of $G = SU(n)$ --- we shall comment 
more on this below.   

Equivariant gerbes were also studied in the paper \cite{BXZ} 
of Behrend, Xu and Zhang;   the authors constructed a bundle gerbe using the path-fibration and 
equipped it with an equivariant connection and curving. 
This construction was followed shortly by a paper of Gawedski and Reis \cite{Gawedzki-Reis2} giving a generalisation  of Meinrenken's  construction to the case of 
non-simply connected groups. 

Finally we come to the case of interest in this paper which is the construction of a local bundle gerbe 
over $SU(n)$.  As mentioned above this example first appeared in the paper \cite{Meinrenken} of Meinrenken 
and was later discussed also by Mickelsson \cite{Mic}. 
In this construction a local bundle gerbe was defined over the open cover 
of $SU(n)$ by open sets $U_z$ consisting of unitary matrices for which $z$ is not 
an eigenvalue.

We will show how to remove the dependence on the local cover in Meinrenken's and Mickelsson's construction 
and to generalise it to any group $G$  which is one of the following: 
the unitary group $U(n)$ (more generally the group $U(H)$ 
of unitary operators on some finite dimensional Hilbert space $H$), the (diagonal) torus  $T = \TT^n \subset U(n)$, or one of the 
Banach Lie groups $U_p(H)$ for $H$ an infinite dimensional separable  complex Hilbert space.  These are defined in more detail below in Section \ref{sec:basic}. 
A common feature of all these groups is that they have natural unitary representations on finite or infinite dimensional Hilbert spaces and so, for convenience, we refer to them 
as unitary groups.  If $G = U(n)$ or $G = U_p(H)$ then $H^3(G,\ZZ) = \ZZ$ and the image in real 
cohomology of the generator of $H^3(G,\ZZ)$ is represented by the basic three-form~\eqref{eq: basic 3-form on SU(n)} 
(note the trace in~\eqref{eq: basic 3-form on SU(n)} makes sense for $g\in U_p(H)$ if $1\leq p\leq 2$).  
Our main result is an explicit construction of a natural connection and curving on this bundle gerbe which simultaneously covers all 
of the above examples of unitary groups.  
Of particular interest is the 
fundamental role played by the holomorphic functional calculus for operators which allows us to obtain explicit, albeit slightly complicated, formulae for the curving (note 
that the functional calculus was also used in \cite{Mic} to construct trivialisations of a certain gerbe).  
We show explicitly that the three-curvature of this connection and curving is $2\pi i$ times the basic three-form on $G$ and hence that  the basic three-form represents the real Dixmier-Doudy class of the basic bundle gerbe. 
 The construction of the basic bundle gerbe that we give is manifestly equivariant but we do not attempt to construct an equivariant connection and curving.   
 
We should also comment on the relation of the basic gerbe to the `index gerbe' considered in a number of papers, beginning with 
\cite{CMM} and further studied in \cite{Bunke, Lott}.  Carey and Wang in \cite{CW} gave a simpler treatment, which 
clarified the relationship of this index gerbe to the families index theorem.   
In \cite{Mic} (see also \cite{CM}) it is explained how the basic gerbe on $G$ can be regarded as a special case of the index gerbe: 
one regards points in $G$ as specifying holonomies of connections on the trivial $G$-bundle over $S^1$; these connections 
can then be coupled to the Dirac operator $-i\frac{d}{dz}$ on $S^1$, giving a family of self adjoint operators parametrised by $G$.  
In the above cited works on the index gerbe it is evident that having explicit formulae for connections and curving is of some importance.   

Finally, in summary, the first two sections give necessary background on bundle gerbes and holomorphic functional calculus. The so-called basic bundle gerbe is introduced
in the next section and this is followed by a construction of a connection and curving on it.  The Weyl map is used to prove that the basic three-form represents the real Dixmier-Douady class of the basic bundle gerbe. In the final section we comment briefly upon the equivariant case.  A number of complicated proofs using the holomorphic functional calculus have been relegated to the appendices.

\section{Background material}
\label{sec: background material}
\subsection{Bundle gerbes}
We review briefly the definition of bundle gerbes \cite{Mur96Bundle-gerbes}. It will be convenient to use hermitian line bundles
in the definition instead of $U(1)$ principal bundles. These two approaches are, of course, equivalent.
Let $\pi \colon Y \to M$ be a surjective submersion and let $Y^{[p]}$ be the $p$-fold fibre product
$$
Y^{[p]} = \{(y_1, \dots, y_p) \mid \pi(y_1) = \dots = \pi(y_p) \} \subset Y^p .
$$
For each $i=1, \dots, p+1$ define the projection  $\pi_i \colon Y^{[p+1]} \to Y^{[p]}$ to be the map that omits the $i$-th element. 
If $J \to Y^{[p]}$ is a line bundle we define a new line bundle $\delta(J) \to Y^{[p+1]}$
by
$$
\delta(J) = \pi_1^*(J) \otimes \pi_2^*(J)^* \otimes \pi_3^*(J) \otimes \cdots.
$$
It is straightforward to check that $\delta(\delta(J)) $ is canonically trivial. Moreover if $\sigma$ is a
section of $J$ then $\delta(\delta(\sigma)) = 1$ under this canonical trivialisation.
We then have:
\begin{definition}
A {\em bundle gerbe} on $M$ is a pair $(L, Y)$ where $\pi \colon Y \to M$ is a surjective submersion and
$L \to Y^{[2]}$  a line bundle such that:
\begin{enumerate}
\item there is a non-zero section $s$ of $\delta(L) \to Y^{[3]}$, and
\item $\delta(s) = 1 $ as  a section of $\delta(\delta(L))$.
\end{enumerate}
\end{definition}

If $L$ is a hermitian line bundle and $s$ has length one we call $(L, Y)$ a
{\em hermitian bundle gerbe} on $M$.
Notice that if $(y_1, y_2, y_3) \in Y^{[3]}$ then a vector 
$$
s(y_1, y_2, y_3) \in L_{(y_2, y_3)} \otimes L_{(y_1, y_3)}^* \otimes L_{(y_1, y_2)}
$$
of length one defines a hermitian isomorphism
\begin{equation}
\label{eq:bgmult}
L_{(y_1, y_2)} \otimes L_{(y_2, y_3)} \to L_{(y_1, y_3)}
\end{equation} 
called the {\em bundle gerbe multiplication}. If $(y_1, y_2, y_3, y_4)  \in Y^{[4]}$ there
are two natural ways to compose the bundle gerbe multiplication and define an isomorphism:
$$
L_{(y_1, y_2)} \otimes L_{(y_2, y_3)} \otimes L_{(y_3, y_4)} \to   L_{(y_1, y_4)}.
$$ 
If these two isomorphisms agree we call the bundle gerbe multiplication {\em associative} and this is 
equivalent to the condition $\delta(s) = 1$. 

Associated to any bundle gerbe is a class in $H^3(M, \ZZ)$ called the {\em Dixmier-Douady class} of the bundle gerbe. We will omit
its derivation, which can be found in \cite{Mur96Bundle-gerbes}, as we do need it in the discussion below. We do, however, need 
to understand its image in real cohomology which can be defined  as follows.

Let $\Omega^q(Y^{[p]})$ be the space of differential $q$-forms on $Y^{[p]}$.  Define a homomorphism 
$$
\delta \colon \Omega^q(Y^{[p]}) \to \Omega^q(Y^{[p+1]}) \quad\text{by}\quad \delta =  \sum_{i=1}^{p+1} (-1)^{i-1} \pi_i^* .
$$

These maps form the {\em fundamental  complex}
$$
0 \to \Omega^q(M) \stackrel{\pi^*}{\to} \Omega^q(Y) \stackrel{\delta}{\to}  \Omega^q(Y^{[2]}) \stackrel{\delta}{\to}
\Omega^q(Y^{[3]})    \stackrel{\delta}{\to} \dots
$$
which is exact \cite{Mur96Bundle-gerbes}.
If $(L, Y)$ is a bundle gerbe on $M$ then a {\em bundle gerbe connection} is a
connection $\nabla$ on $L$ such that $s$ is covariantly constant
for $\delta(\nabla)$.  Equivalently the connection $\nabla$ commutes with the bundle gerbe multiplication 
\eqref{eq:bgmult}. If $\nabla$ is a bundle gerbe connection and $F_\nabla$ is
its curvature then $\delta(F_\nabla) = 0$ so that $F_\nabla = \delta(f)$ for some two-form $f \in \Omega^2(Y)$. 
A choice of such an $f$ is
called a {\em curving} for $\nabla$. From the exactness of the fundamental complex we see that the 
curving is only unique up to addition of two-forms pulled back to $Y$ from $M$.  
Given a choice of curving $f$ we have  $\delta(df) = d \delta(f) =
dF_\nabla = 0$ so that $df = \pi^*(\omega)$ for a closed three-form $\omega$ on $M$ called the \emph{three-curvature}
of $\nabla$ and $f$.  The de Rham class
$$
\left[ \frac{1}{2\pi i} \omega \right] \in H^3(M, \RR)
$$
is an integral class which is the image in real cohomology of the Dixmier-Douady class of  $(L, Y)$. For convenience
let us call this the real Dixmier-Douady class of $(L, Y)$.

\subsection{Holomorphic functional calculus}
\label{sec:hfc}

We briefly 
recall the main features of the holomorphic functional calculus and refer the reader to \cite{DunSch} for more details. 
Let $B(H)$ denote the Banach algebra of bounded operators on a Hilbert space $H$.  Suppose that $T\in B(H)$ and that the spectrum 
$\spec(T)$ of $T$ is a compact subset of the complex plane.  Given  
a complex valued function $f$, holomorphic 
on an open neighbourhood $U$ of $\spec(T)$, we can 
define a new operator $f(T)$ by the contour integral 
$$ 
f(T) = \frac{1}{2\pi i}\oint_C f(\xi)(\xi - T)^{-1}d\xi. 
$$ 
Here $C$ is a contour 
surrounding $\spec(T)$ which lies entirely in $U$ and is always taken to be oriented 
counter-clockwise. Because the resolvent $(\xi - T)^{-1}$ is holomorphic away from $\spec(T)$ this 
definition is independent of continuous deformations of $C$.   It is an easy consequence of Cauchy's theorem 
that if $g$ 
is another complex valued function, holomorphic on the same 
open neighbourhood $U$, then $(fg)(T) = f(T)g(T)$.  

As an example of this formula suppose that $\lambda$ is 
an isolated point of $\spec(T)$ which  
is an eigenvalue of $T$.  Then the orthogonal projection 
$P_{\lambda}$ onto the $\lambda$ eigenspace of $T$, $E_{(T, \lambda)}$, 
is given by the contour integral 
\begin{equation} 
\label{eq: contour integral} 
P = \frac{1}{2\pi i}\oint_{C_\lambda} (\xi - T)^{-1}d\xi 
\end{equation}
where $C_{\lambda}$ is a small circle centred at $\lambda$ 
which contains no other point of $\spec(T)$.  It is instructive to consider this for the case when $T = g$ is 
a unitary matrix.   
Since we can write $g$ as a sum $\sum \lambda_i P_i$ 
where $P_i$ denotes the orthogonal projection onto the 
$\lambda_i$-eigenspace we see that the resolvent 
$(\xi - g)^{-1}$ can be written as a sum 
$$
(\xi - g)^{-1} = \sum (\xi - \lambda_i)^{-1}P_i.
$$ 
This gives an effective way to evaluate the contour 
integral~\eqref{eq: contour integral} since this can now 
be written as  
$$ 
\frac{1}{2\pi i}\oint_{C_{\lambda}} \sum (\xi - \lambda_i)^{-1}P_id\xi 
$$ 
and all one has to do is evaluate the 
contour integrals $\frac{1}{2\pi i}\oint_{C_{\lambda}} (\xi - \lambda_i)^{-1}d\xi$ 
which can easily be done.

In a similar fashion assume that $C$ is a contour encircling eigenvalues 
$\lambda_1, \dots, \lambda_r$ of $T$. Then 
\begin{equation} 
\label{eq:proj} 
P = \frac{1}{2\pi i}\oint_{C} (\xi - T)^{-1}d\xi 
\end{equation}
gives the projection  onto the sum of the eigenspaces:
$$
E_{(T, \lambda_1)} \oplus \cdots \oplus E_{(T, \lambda_r)}
$$


\section{The basic bundle gerbe}
\label{sec:basic}
We now give a global version of the construction in  \cite{Meinrenken, Mic} of a bundle gerbe 
on a group of unitary operators.  Suppose that $G$ is one of the following groups: 
the unitary group $U(n)$ (more generally the group $U(H)$ 
of unitary operators on some finite dimensional Hilbert space $H$), the (diagonal) torus  $T = \TT^n \subset U(n)$, or one of the 
Banach Lie groups $U_p(H)$ for $H$ an infinite dimensional separable 
complex Hilbert space.  Recall \cite{Quillen} that for $p\geq 1$ the groups $U_p(H)$ are defined 
to be the groups of unitary operators on $H$ differing from the identity by an 
operator in the Schatten ideal $\cL_p(H)$.  For more details on the ideals 
$\cL_p(H)$ we refer the reader to \cite{Simon}.  Notice that in all cases elements of $G$ act
as unitary operators on a  Hilbert space $H$.   For convenience of presentation we shall
refer to the space that $G$ acts on as $H$ even when it might be more natural to call it $\CC^n$. 

Denote the identity in $G$ by $1$ and
note that if $g \in G$ then its spectrum $\spec(g)$ (when we consider $g$ as an operator on $H$) is a subset of the circle $U(1) \subset \CC$.  
If $H$ is finite dimensional then the spectrum of $g \in G$   is of course finite, and hence discrete. 
Suppose that $H$ is infinite dimensional so that  $g\in U_p(H)$ for some $p \geq 1$.    Write $g = 1 + A$ where 
$A$ belongs to the appropriate ideal $\cL_p(H)$.  If $\lambda \in \spec(g)$ then 
$g - \lambda$ is not invertible and hence $A - (\lambda - 1)$ is not invertible.  This 
means $\lambda - 1$ belongs to the spectrum of the compact operator $A$.   
The spectrum $\spec(A)$ of the compact operator $A$ is a set with no non-zero accumulation
points, and if $\mu\in \spec(A)$ is non-zero then $\mu$ is an eigenvalue of 
finite multiplicity.  It follows that the spectrum $\spec(g)$ of $g$ is a subset of 
$U(1)$ which has at most one accumulation point at $1\in U(1)$.  All elements of 
the spectrum of $g$ not equal to $1$ are eigenvalues of finite multiplicity.    
Denote by $U_0(1)$ the set $U(1)$ with the identity element, $1$, removed.

Define
$$
Y = \{ (z, g) \mid z  \notin \spec(g) \cup \{ 1 \}  \} \subset  U_0(1) \times G
$$
and let $\pi \colon Y \to G$ be the projection $\pi(z, g) = g$.
We will identify $Y^{[p]}$ with the subset of $ U_0(1)^p \times G$ determined by
$$
Y^{[p]}  = \{ (z_1, \dots, z_p, g) \mid z_1,\dots, z_p   \notin \spec(g) \cup \{ 1 \}  \} \subset  U_0(1)^p \times G
$$
We put an ordering on $U_0(1)$  as follows.   If $z_1, z_2   \in U_0(1)$ we say that $z_1 > z_2$ if we can positively rotate $z_2$ into $z_1$ without passing through $1$.   For any pair $z_1, z_2 \in U_0(1)$ we say that  $\lambda$ is {\em between} $z_1$ and $z_2$ if $\lambda$ is in the
connected component of $U(1) - \{z_1, z_2\}$ not containing $1$.  
  Call a point $(z_1, z_2, g) \in Y^{[2]}$ {\em positive} if there are eigenvalues of $g$ between $z_1$ and $z_2$ and $z_1 > z_2$, 
{\em null} if there are no eigenvalues of $g$ between $z_1$ and $z_2$ and {\em negative} if 
there are eigenvalues of $g$ between $z_1$ and $z_2$ and $z_1 < z_2$. Note that if $(z_1, z_2, g)$ is
positive (negative) then $(z_2, z_1, g)$ is negative (positive).  Denote the corresponding 
subsets of $Y^{[2]}$ by $Y^{[2]}_+$, $Y^{[2]}_0$ and $Y^{[2]}_-$.  

Let  $E_{(g, \lambda)}$ be the $\lambda$ eigenspace of $g$ for $\lambda \in U_0(1)$. 
If  $(z_1, z_2, g)$ is in $Y^{[2]}_+$ we define
\begin{equation}
\label{eq: direct sum of e-spaces}
E_{(z_1, z_2, g)} = \bigoplus_{z_1 >\lambda > z_2} E_{(g,\lambda)}.
\end{equation}
By construction,
$E_{(z_1, z_2, g)}$ has finite dimension so we can define the top exterior power
$$
L_{(z_1, z_2, g)}  = \det(E_{(z_1, z_2, g)})
$$
If $(z_1, z_2, g) \in Y^{[2]}_0$ we define
$$
L_{(z_1, z_2, g)}  = \CC
$$
and if $(z_1, z_2, g) \in Y^{[2]}_-$ we define
$$
L_{(z_1, z_2, g)}  = \det(E_{(z_2, z_1, g)})^*.
$$

We want to show that $L \to Y^{[2]}$ is a smooth locally trivial hermitian line bundle. It suffices to show 
that $E \to Y^{[2]}_+$ is smooth and locally trivial and by standard results this follows if we can show that 
$P\colon Y^{[2]}_+\to B(H)$ is smooth where $P$ is the orthogonal projection onto $E$. 
To  show this last fact we will use the 
holomorphic functional calculus from Section \ref{sec:hfc}. 

\begin{center}
\begin{figure}
\includegraphics[scale=0.4]{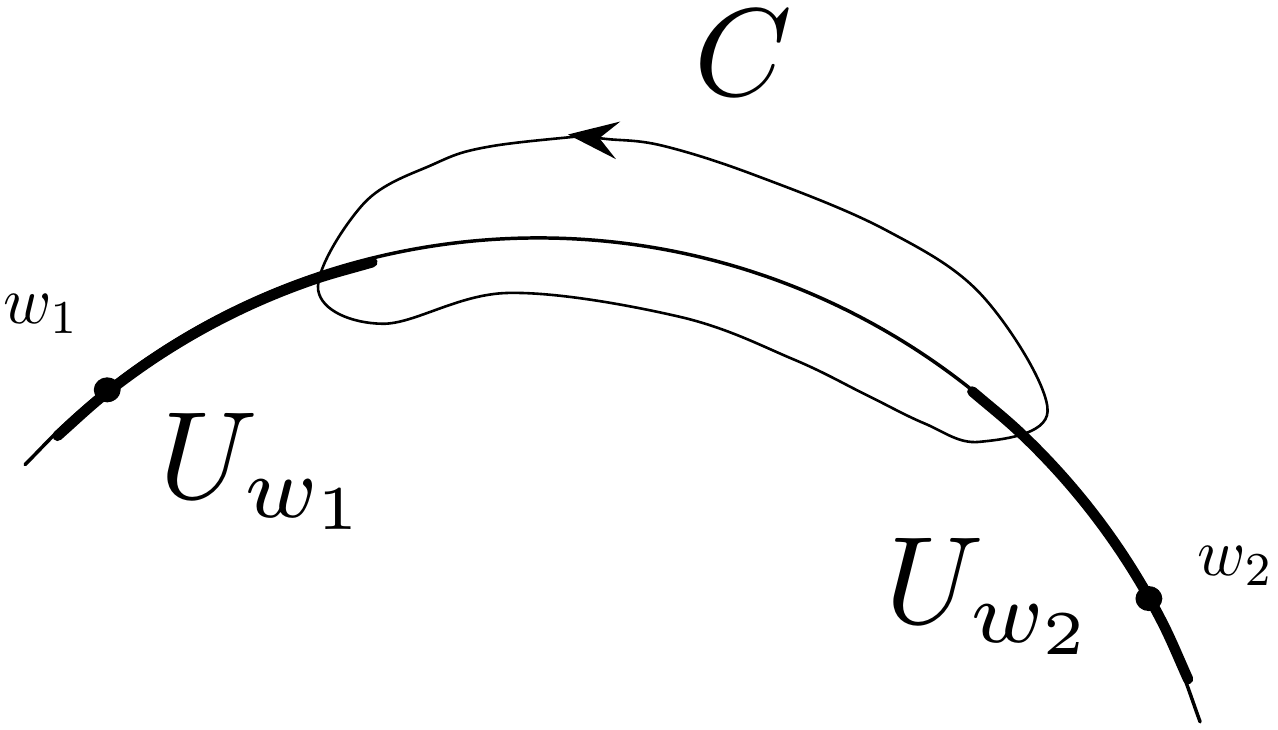}
\caption{The contour $C$.}
\label{fig:contour}
\end{figure}
\end{center}

 Consider the continuous map
\begin{equation}
\begin{array}{ccc}
U_0(1) \times U_0(1) \times G & \to & B(H) \times B(H)\\
(z_1, z_2, g) & \mapsto & (z_1 - g, z_2 - g) \\
\end{array}
\end{equation}
As $Y^{[2]}$ is the pre-image under this of the open set of 
pairs of invertible operators  in $B(H)\times B(H)$ it follows that $Y^{[2]}$ is open in $U_0(1) \times U_0(1) \times G$. If $(w_1, w_2, h) 
\in Y^{[2]}_+$ we can therefore find connected open sets $U_{w_1}, U_{w_2}$ and $U_{h}$ containing 
$w_1$, $w_2$ and $h$ respectively, with the property that if 
$$
(z_1, z_2, g) \in U_{w_1} \times U_{w_2} \times U_{h} \subset Y^{[2]}
$$
then neither of $z_1$ or $z_2$ is in the spectrum of $g$. In particular, if $g \in U_{h}$  and $C$ 
is any  anti-clockwise contour cutting $U(1)$ once in $U_{w_1}$ and once in $U_{w_2}$ 
then $C$ must encircle all the eigenvalues of $g$ between $z_1$ and $z_2$. Fix such a 
contour as in Figure \ref{fig:contour}.
Let 
\begin{equation}
\label{eq: open nbhood}
U_{(w_1, w_2, h)} = U_{w_1} \times U_{w_2} \times U_{h}.
\end{equation}
The projection map restricted to $U_{(w_1, w_2, h)}$
$$
P \colon U_{(w_1, w_2, h)} \to B(H)
$$
is therefore  given by the contour integral formula 
\begin{equation} 
\label{eq: reg expression for P} 
P(z_1, z_2, g) = \frac{1}{2\pi i}\oint_{C} (\xi - g)^{-1}d\xi .
\end{equation}
As we can differentiate through the contour integral and the integrand is 
smooth it follows that $P$ is smooth on $U_{(w_1, w_2, h)}$ and hence on all 
of $Y^{[2]}$.  Moreover 
$
\tr(P) = \dim(E_{(z_1, z_2, g)})
$
must be constant and equal to $\dim(E_{(w_1, w_2, h)})$ so that 
$
U_{(w_1, w_2, h)} 
$
must lie entirely in one of $Y^{[2]}_+$, $Y^{[2]}_0$ or $Y^{[2]}_-$ and hence they 
must all be open.  It follows that $P \colon Y^{[2]}_+ \to B(H)$ is smooth and we have shown that:
\begin{proposition}
$L$  is a smooth and locally trivial line bundle on $Y^{[2]}$.
\end{proposition}

We want to now prove that $L$ is a bundle gerbe. The basic fact that this uses is that if $ V \oplus W $ is 
a direct sum of vector spaces then wedge product defines  a canonical isomorphism
$$
\det(V) \otimes \det(W) = \det(V \oplus W)
$$
However the  proof is made complicated by the various cases that arise. To handle some of these we make the following definitions. 
Let $\Sigma_r$ denote the symmetric group on $r$ letters. This   acts naturally on $Y^{[r]}$ by permuting the $y_1, \dots, y_r$
in $(y_1, \dots, y_r, g)$.   If $J \to Y^{[r]}$ is a 
line bundle we say that $J$ is  \emph{antisymmetric}  if for all $\rho \in \Sigma_r$ we have isomorphisms
$$
J \simeq \rho^*(J)^{\sgn(\rho)}
$$
where we denote by $\rho \colon Y^{[r]} \to Y^{[r]}$ the map induced by the permutation $\rho$. Also  if $V$ is a vector space, we are using the notation $V^{-1}$ for $V^*$. 
If $0 \neq v \in V$ a one-dimensional vector space define $v^* \in V^*$ by $v(v^*) = 1$ and also denote
$v^*$ by $v^{-1}$.  If $J$ is antisymmetric and $\psi$ is a non-vanishing section of $J$ then it makes sense to
define $\psi$ to be antisymmetric if $\psi = \rho^*(\psi)^{\sgn(\rho)}$ for all $\rho \in \Sigma_r$.
\begin{lemma} 
\label{lemma:antisymmetric}
If $J \to Y^{[r]}$ is  antisymmetric then $\delta(J)$ is  antisymmetric. In such a  case if $\psi $ is an antisymmetric
section of $J$ then $\delta(\psi)$ is an antisymmetric section of $\delta(J)$. 
\end{lemma}
\begin{proof} 
We have
$$
\delta(J)_{(x_1, x_2, \dots, x_{r+1})} = J_{(x_2, x_3, \dots, x_{r+1})} \otimes J_{(x_1, x_3, \dots, x_{r+1})}^* \otimes 
J_{(x_1, x_2, x_4, \dots  \dots , x_{r+1})}\otimes \cdots 
$$
Inspection shows that any side by side transposition changes the right hand side to its dual so the result follows.
The same method works for $\psi$. 
\end{proof}
\noindent Notice that, by construction, $L \to Y^{[2]}$ is an antisymmetric hermitian line bundle. 
Consider $(z_1, z_2, z_3, g) \in Y^{[3]}$. There is a  $\rho \in \Sigma_3$ such that 
$z_{\rho(1)} \geq z_{\rho(2)} \geq z_{\rho(3)}$. Define $(z_1, z_2, z_3, g)$ to be of type 
$(1, 1)$ if there are eigenvalues of $g$ between $z_{\rho(1)} $ and $ z_{\rho(2)}$ and 
also between $z_{\rho(2)}$ and $z_{\rho(3)}$, type $(1, 0)$ if there are eigenvalues of $g$  between $z_{\rho(1)} $ 
and $ z_{\rho(2)}$ but not  between $z_{\rho(2)}$ and $z_{\rho(3)}$, type $(0, 1)$ if there are no 
eigenvalues of $g$  between $z_{\rho(1)} $ and $ z_{\rho(2)}$ but there are some  between 
$z_{\rho(2)}$ and $z_{\rho(3)}$ and type $(0,0)$ if there are no eigenvalues of $g$ between 
$z_{\rho(1)} $ and $ z_{\rho(3)}$. Notice that although $\rho$ is not unique if some of the $z_i$ are equal, nevertheless these 
definitions still make sense.  Denote by $Y^{[3]}_{(i,j)}$ the subset of $Y^{[3]}$ consisting
of elements of type $(i, j)$ for each $i, j = 0, 1$. Each $Y^{[3]}_{(i,j)}$ is invariant under $\Sigma_3$ and  
is a union of connected components of  $Y^{[3]}$ and hence open.

Consider  $(z_1, z_2, z_3, g) \in Y^{[3]}_{(1, 1)}$ with $z_1 > z_2 > z_3$. 
We have 
\begin{equation}
\left(  \bigoplus_{ z_1 > \lambda > z_2} E_{(g, \lambda)} \right) \oplus \left(  \bigoplus_{ z_2 > \lambda > z_3} E_{(g, \lambda)} \right) = \left(  \bigoplus_{ z_1 > \lambda > z_3} E_{(g, \lambda)} \right) 
\end{equation}
so that
\begin{equation}
\label{eq:directsum}
E_{(z_1, z_2, g)} \oplus E_{(z_2, z_3, g)} =  E_{(z_1, z_3, g)} 
\end{equation}
and hence wedge product gives an isomorphism
\begin{equation}
\label{eq:mult}
L_{(z_1, z_2, g)} \otimes L_{(z_2, z_3, g)}
=  L_{(z_1, z_3, g)}.
\end{equation}
Moreover this  defines a smooth section
\begin{equation}
\label{eq:bgsec}
s(z_1, z_2, z_3, g) \in \delta(L)_{(z_1, z_2, z_3, g)} = 
L_{(z_2, z_3, g)} \otimes L^*_{(z_1, z_3, g)} \otimes L_{(z_1, z_2, g)} = \CC
\end{equation}
at points satisfying $z_1 > z_2 > z_3$.  We can extend this to a section $s$ of $\delta(L)$ over all of $Y^{[3]}_{(1, 1)}$
by requiring antisymmetry. 

Consider  $(z_1, z_2, z_3, g) \in Y^{[3]}_{(1, 0)}$ with $z_1 > z_2 \geq z_3$.  Then by definition
\begin{equation}
\label{eq:cond1}
L_{(z_2, z_3, g)} = \CC
\end{equation}
and moreover $E_{(z_1, z_2, g)} = E_{(z_1, z_3,g)}$ so that
\begin{equation}
\label{eq:cond2}
L_{(z_1, z_2, g)} = L_{(z_1, z_3, g)}.
\end{equation}
Thus we have 
$$
\delta(L)_{(z_1, z_2, z_3, g)} =   L_{(z_2, z_3, g)} \otimes L^*_{(z_1, z_3, g)} \otimes L_{(z_1, z_2, g)} 
= \CC \otimes  L^*_{(z_1, z_3, g)} \otimes L_{(z_1, z_3, g)} = \CC
$$
so we can define
$$
s(z_1, z_2, z_3, g) \in \delta(L)_{(z_1, z_2, z_3, g)}.
$$
Notice that on the connected component containing $(z_1, z_2, z_3, g) $ the conditions in equations \eqref{eq:cond1} and 
\eqref{eq:cond2} will be satisfied for all points so that $\psi$ extends to a smooth section on all of that connected component. 
A similar argument can be applied at the other points of $Y^{[3]}_{(1, 0)}$ to give an antisymmetric section $s$.  
Clearly
the same type of argument can be applied in the case of $Y^{[3]}_{(0,1)}$. 
On the remaining subset $Y^{[3]}_{(0, 0)}$ all the line bundles are trivial so there is an obvious section. 
We conclude that there is a  naturally defined antisymmetric section $s$ of $\delta(L)$ over all of  $Y^{[3]}$.

We have already remarked that   $\delta(s) =1 $ is equivalent to the bundle gerbe multiplication 
\eqref{eq:bgmult} being associative.  In the case that $z_1 > z_2 > z_3 > z_4$ and there are eigenvalues of $g$ between each of the 
consecutive $z$'s this follows because the wedge product
is associative.  As $s$ and $L$ are both antisymmetric it follows from Lemma \ref{lemma:antisymmetric}
that $\delta(s)$ is antisymmetric, so that $\delta(s)(z_1, z_2, z_3, z_4, g)  = 1$ whenever 
there are  eigenvalues of $g$ between each  consecutive $z$'s, regardless of their ordering.  The other cases can be dealt with in a similar
fashion.  
We conclude that $L \to Y^{[2]}$ is a hermitian bundle gerbe on $G$.

\section{A bundle gerbe connection and its curvature}

The basic bundle gerbe has a canonical bundle gerbe
connection constructed as follows (this construction is mentioned in \cite{CMM} and 
discussed also in \cite{Gawedzki-Reis}). Over $Y^{[2]}_+$ we have
$$
E \subset H \times Y^{[2]}_+
$$
and we can use the  orthogonal
projection
$$
P \colon H \to E
$$
to project the trivial connection on $H \times Y^{[2]}_+
$ to a  connection $\nabla$ on $E$.  
Over $Y^{[2]}_+$ this defines a connection 
$\det(\nabla) $ on $L$, over $Y^{[2]}_-$ we take the corresponding dual connection and over
$Y^{[2]}_0$ we take the flat connection. Notice that as $L$ is antisymmetric we can 
consider the corresponding notion of antisymmetry for a connection and by construction $\det(\nabla)$ is  antisymmetric. 

We wish to show this is a bundle gerbe connection. Consider a connected component $X$
of $Y_{(1,1)}^{[3]}$ containing some $(z_1, z_2, z_3, g)$ for which  $z_1 > z_2 > z_3$. From equation \eqref{eq:directsum}  we have the orthogonal direct sum
\begin{equation}
E_{(z_1, z_3, g)} = E_{(z_1, z_2, g)} \oplus E_{(z_2, z_3, g)} \subset H.
\end{equation}
It follows that over $X$ we have
\begin{equation}
\label{eq:Esplit}
\pi^*_2(E) = \pi^*_3(E) \oplus \pi^*_1(E) \subset Y^{[3]} \times H
\end{equation}
and moreover that we have
$$
\pi^*_2(\nabla) = \pi^*_3(\nabla) \oplus \pi^*_1(\nabla)
$$
because equation \eqref{eq:Esplit} is an orthogonal splitting and the 
connections are defined by orthogonal projections.   Hence over $X$ the connection respects the
wedge product isomorphism
$$
\pi^*_2(L) =  \pi^*_3(L) \otimes \pi^*_1(L)
$$
and hence the section $s$ is covariantly constant for the connection on $L$ over $X$. By antisymmetry
the section $s$ is covariantly constant for the connection over all of $Y^{[3]}_{(1, 1)}$. 
The other parts of $Y^{[3]}$ can be dealt with in a similar fashion. We conclude that the connection
we have constructed is a bundle gerbe connection.

Over $Y^{[2]}_+$ the  curvature two-form $F_{\det(\nabla)}$ of the bundle gerbe connection
will be equal to the trace $\tr(F_\nabla)$
of the curvature $F_\nabla$ of the connection
$\nabla$ on the vector bundle $E$.  It is
a simple calculation to see that this can be computed in terms of the projections
$P$ as
\begin{equation} 
\label{eq: curv expression}
F_{\det(\nabla)} = \tr(PdPdP). 
\end{equation}
By the antisymmetry of $\det(\nabla)$ the sign of $F_{\det(\nabla)}$ will change on $Y^{[2]}_-$ and it 
will vanish on $Y^{[2]}_0$.

We now  explain how the holomorphic functional calculus  from Section \ref{sec:hfc}
can be used to derive a contour integral expression for the 
curvature two-form $F_{\det(\nabla)}$ described above.  This will be useful in defining a curving
for $F_{\det(\nabla)}$.
We can use the local expressions~\eqref{eq: reg expression for P}   
for the projections $P$ to write down a contour integral formula 
for $F_{\det(\nabla)}$; a priori this will be a triple contour integral, 
however it is possible to obtain a simpler formula as in the following 
Proposition.  
Let $G$ be one of the groups $U(H)$ for $H$ a finite dimensional 
complex Hilbert space, or $U_p(H)$ for $H$ an infinite dimensional 
complex Hilbert space and $1\leq p\leq 2$.   For any 
$(z_1,z_2,g)\in Y^{[2]}_+$ choose a closed contour $C_{(z_1,z_2,g)}$ 
enclosing all the eigenvalues of $g$ between $z_1$ and $z_2$, oriented 
counter clockwise.  Then we have 
\begin{proposition} 
\label{prop: expression for bg curvature} 
The restriction of the curvature $F_{\det(\nabla)}$ to 
$Y^{[2]}_+$ is given by
\begin{equation} 
\label{eq: contour integral formula for curvature} 
F_{\det(\nabla)}(z_1,z_2,g) = \frac{1}{4\pi i}\oint_{C_{(z_1,z_2,g)}}\tr\left((\xi - g)^{-1}dg(\xi - g)^{-2}dg\right) d\xi 
\end{equation}
\end{proposition} 
Notice that on the right hand side of  this formula we have committed an abuse of notation and denoted by $dg$ the derivative at $(z_1,z_2,g)$ 
of the projection map $Y^{[2]}\to G$.  
Notice also  that the curvature two-form has no components in the $U(1)$ directions.  
We also need to comment on this expression in the infinite dimensional case when $G = U_p(H)$.  
If $1\leq p\leq 2$ then the trace in~\eqref{eq: contour integral formula for curvature} makes 
sense, since $(\xi - g)^{-1}dg(\xi - g)^{-2}dg$ is a product of one-forms taking 
values in $\cL_p(H)$, and hence is a two-form taking values in the trace class 
operators on $H$.  
The proof of Proposition~\ref{prop: expression for bg curvature} 
is deferred to Appendix~\ref{sec: first proof}.

\section{A curving} 

We would now like to find a curving for this bundle gerbe connection, 
in other words we would like to find a two-form $f$ on $Y$ such that 
$$ 
F_{\det(\nabla)} = \delta(f)  = \pi_1^*f - \pi_2^*f 
$$ 
on $Y^{[2]}$.  We can do this as follows.  For each complex 
number $z$ with $|z| = 1$ let $R_z$ denote the closed ray from the origin through $z$. 
For any $z \in U_0(1)$  we define a branch of the logarithm, $\log_z \colon \CC - R_z \to \CC $,  
by making the cut along $R_z$ and also setting $\log_z(1) = 0$.  If 
$z_1 > z_2$ write $(z_1,z_2)$ for the set of $\xi\in U(1)$ between 
$z_1$ and $z_2$ in the sense of Section~\ref{sec:basic}.   
It is easy to check that if 
$z_1> z_2$ then 
\begin{equation} 
\label{eq: property of log_z} 
\log_{z_1} \xi - \log_{z_2} \xi =\begin{cases} 
2\pi i & \text{if}\ \xi/|\xi| \in (z_1,z_2) \\ 
0 & \text{otherwise.} 
\end{cases} 
\end{equation}
For any $(z,g)\in Y$, choose an anti-clockwise 
closed contour $C_{(z,g)}$ in $\CC - R_z$ which encloses $\spec(g)$. 

\begin{theorem} 
\label{thm: main thm}
Suppose that $G = U(H)$ for some finite dimensional complex 
Hilbert space $H$.   
Define a two-form $f$ on $Y$ by 
\begin{equation} 
\label{eq: curving formula}
f(z,g) = \frac{1}{8\pi^2}\oint_{C_{(z,g)}}
\log_z \xi \tr\left((\xi - g)^{-1}dg(\xi - g)^{-2}dg\right) d\xi.
\end{equation}
Then we have:

\noindent(a) the two-form $f$ is a curving for the bundle gerbe 
connection $\det(\nabla)$ in the sense that it 
satisfies $\delta(f) = F_{\det(\nabla)}$ and, \\
\noindent(b)
the 3-curvature of the bundle gerbe connection $\det(\nabla)$ 
and curving $f$ is the  three-form 
$$ 
-\frac{i}{12\pi}\tr(g^{-1}dg)^3. 
$$ 
\noindent{(c)} The real Dixmier-Douady class of the basic bundle gerbe $(L, Y)$ is 
$$ 
-\frac{1}{24\pi^2}\tr(g^{-1}dg)^3. 
$$
\end{theorem} 
Again we have abused notation on the right hand side of \eqref{eq: curving formula} and denoted by $dg$ what is really the derivative at $(z, g)$ of the projection from $Y^{[2]}$ to $G$.


We need to show first of all that~\eqref{eq: curving formula} 
defines a smooth two-form on $Y$.  Observe that $f(z,g)$ is independent of 
the choice of contour $C_{(z,g)}$, since the integrand is holomorphic in $\xi$.  Fix $(w, h)\in Y$.  We can choose 
an open neighbourhood  $U_{w}$ of $w$ in $U_0(1)$, and an open neighbourhood $U_{h}$ of $h$ in  
$G$ such that $U_{(w, h)} = U_{w}\times U_{h}$ is an open neighbourhood of $(w, h)$ in $Y$.  After perhaps shrinking $U_w$ a little 
we conclude that if $g \in U_h$ then $\spec(g)$ does not intersect $\overline U_w$. 
Therefore we can  find a contour $C$ which contains the spectrum of every $g \in U_h$ and 
lies inside
$$
\CC - \overline{\bigcup_{z \in U_w} R_z}.
$$
Such a  $C$ satisfies the requirements to be a $C_{(z, g)}$ for any $(z, g) \in U_{(w, h)}$.  Moreover for 
any $z \in U_w$ we have $\log_z = \log_w $ on $C$ by~\eqref{eq: property of log_z}. Hence it follows that the 
restriction of $f$ to $U_{(w, h)}$ satisfies 
$$ 
f(z,g) = \frac{1}{8\pi^2}\oint_{C}\log_{w}\xi\tr((\xi - g)^{-1}dg(\xi - g)^{-2}dg)d\xi. 
$$ 
and is clearly smooth on this open set and hence on all of $Y^{[2]}_+$. 

Part (a) of Theorem~\ref{thm: main thm} is  
not too difficult to establish.  It is enough to prove this on the 
open subset $Y^{[2]}_+$ of $Y^{[2]}$.  Recall that on this set $z_1 > z_2$.  
Since the expression~\eqref{eq: curving formula} for $f(z_1,g)$ is independent of the choice of 
contour $C_{(z_1,g)}$ we can replace the contour $C_{(z_1,g)}$ 
with the union 
$$
C_{(z_1,z_2,g)}\cup \tilde{C}_{(z_1,z_2,g)}
$$ 
where the contour $C_{(z_1,z_2,g)}$ is the one described in 
Proposition~\ref{prop: expression for bg curvature} above and the 
contour $\tilde{C}_{(z_1,z_2,g)}$ is a contour surrounding the part of 
the spectrum of $g$ lying in the set $U_0(1)\setminus (z_1,z_2)$.  
Then we have that $\pi_2^*f$ is given by 
\begin{multline*} 
\frac{1}{8\pi^2}\oint_{C_{(z_1,g)}}\log_{z_1}\xi \tr\left( (\xi -g)^{-1} 
dg(\xi - g)^{-2}dg\right) d\xi \\ 
= 
\frac{1}{8\pi^2}\oint_{C_{(z_1,z_2,g)}}\log_{z_1}\xi \tr\left( (\xi -g)^{-1} 
dg(\xi - g)^{-2}dg\right) d\xi \\ 
+ \frac{1}{8\pi^2}\oint_{\tilde{C}_{(z_1,z_2,g)}}\log_{z_1}\xi \tr\left( (\xi -g)^{-1} 
dg(\xi - g)^{-2}dg\right) d\xi.  
\end{multline*}
Likewise we can write the contour $C_{(z_2,g)}$ as a union 
$C_{(z_1,z_2,g)}\cup \tilde{C}_{(z_1,z_2,g)}$ and obtain an expression 
for $\pi_1^*f$ as a sum of contour integrals as above.  It then follows, 
using equation~\eqref{eq: property of log_z} that $\delta(f) = F_{\det(\nabla)}$ 
on $Y^{[2]}_+$.     

The proof of part (b) of Theorem~\ref{thm: main thm}  requires a little preparation, 
in particular we need to make use of some properties of the 
so-called `Weyl map'.  We discuss this in the next section and complete the remainder of the proof   
of Theorem~\ref{thm: main thm} in the Appendix.

It is possible to generalise Theorem~\ref{thm: main thm} to 
include the Banach Lie groups $U_p(H)$ for $1\leq p\leq 2$.  
More precisely we have the following result: 
\begin{theorem} 
\label{thm: mainer thm} 
Let $G$ be one of the Banach Lie groups $U_p(H)$ for $1\leq p\leq 2$.  
Define a two-form $f$ on $Y$ by
\begin{equation} 
\label{eq: formula for curving} 
f(z,g) = \frac{1}{8\pi^2}\oint_{C_{(z,g)}}
\log_z \xi \tr\left((\xi - g)^{-1}dg(\xi - g)^{-2}dg\right) d\xi.
\end{equation}
Then we have:

\noindent(a) the two-form $f$ is a curving for the bundle gerbe 
connection $\det(\nabla)$ in the sense that it 
satisfies $\delta(f) = F_{\det(\nabla)}$ and, \\
\noindent(b)
the 3-curvature of the bundle gerbe connection $\det(\nabla)$ 
and curving $f$ is the  three-form 
$$ 
-\frac{i}{12\pi}\tr(g^{-1}dg)^3. 
$$ 
\noindent{(c)} The real Dixmier-Douady class of the basic bundle gerbe $(L, Y)$ is 
$$ 
-\frac{1}{24\pi^2}\tr(g^{-1}dg)^3. 
$$
\end{theorem} 
The comment about traces following Proposition~\ref{prop: expression for bg curvature} 
applies verbatim in this situation to show that the trace 
in the expression~\eqref{eq: formula for curving} is well-defined.  
The proof that $f$ is a smooth two-form on $Y$ and that $\delta(f) = F_{\det(\nabla)}$ goes through in exactly the same 
manner as above.  We now need to identify the associated 
3-curvature.  

Let $e_1$, $e_2$, $e_3$, \ldots be an orthonormal basis of $H$ 
and let $H_n$ be the subspace of $H$ spanned by $e_1$, $e_2$, 
\ldots, $e_n$.  The algebra $\mathit{gl}(H_n)$ of linear transformations 
of $H_n$ can be identified with a subalgebra of $B(H)$ by sending 
a linear map $T$ to the bounded operator $P_nTP_n$ (here 
$P_n$ denotes the orthogonal projection onto $H_n$).  As pointed 
out in \cite{Quillen}, for any $p,n\geq 1$ there is an inclusion 
$$ 
U(H_n)\subset U_p(H),\ g\mapsto (P_ngP_n - P_n) + 1.  
$$ 
It is clear that the basic gerbe over $U_p(H)$ defined above 
restricts to the basic gerbe over $U(H_n)$ for any $n$.  Also note 
that forms on $U_p(H)$ restrict to forms on $U(H_n)$.  In \cite{Quillen} 
Quillen proves the following result: 
\begin{theorem}[\cite{Quillen} Proposition~7.16] 
\label{thm: Quillens thm}
Let $\omega$ be a form defined on $U_p(H)$ for any $p\geq 1$.  If 
$\omega$ vanishes when restricted to $U(H_n)$ for any $n$, then 
$\omega$ vanishes on $U_p(H)$.   
\end{theorem} 

We can use Theorem~\ref{thm: Quillens thm} to calculate the 3-curvature $\omega$ of the basic 
bundle gerbe with connection $\det(\nabla)$ and curving $f$ on 
$U_p(H)$ for $1\leq p\leq 2$.  Consider the three-form 
$$ 
\omega + \frac{i}{12\pi}\tr(g^{-1}dg)^3 
$$ 
defined on $U_p(H)$.  Under the inclusion $U(H_n)\subset U_p(H)$ 
the left Maurer-Cartan 1-form $\Theta_L = g^{-1}dg$ on $U_p(H)$ pulls back to 
a 1-form on $U(H_n)$ taking values in $\cL_p(H)$.  It is 
easy to see that this pullback of $\Theta_L$ is given by 
the 1-form
$$ 
P_n\theta_LP_n
$$ on $U(H_n)$ where $\theta_L$ denotes the left Maurer-Cartan 
1-form $\theta_L = g^{-1}dg$ on $U(H_n)$.    Since 
$\theta_L$ commutes with $P_n$ it follows that 
 the three-form $-\frac{i}{12\pi}\tr(g^{-1}dg)^3$ 
on $U_p(H)$ restricts to the corresponding three-form $-\frac{i}{12\pi}\tr(g^{-1}dg)^3$ on every 
$U(H_n)$ and so the difference $\omega + \frac{i}{12\pi}\tr(g^{-1}dg)^3$ 
vanishes on every $U(H_n)$.  By Quillen's result this means that 
$$ 
\omega = -\frac{i}{12\pi}\tr(g^{-1}dg)^3. 
$$

\section{The Weyl map $p\colon G/T\times T\to G$} 
\label{sec: Weyl map} 

Suppose that $G$ is a compact, connected Lie group 
and that $T$ is a maximal torus of $G$.  There is a 
canonical map 
$$ 
p\colon G/T\times T\to G 
$$ 
defined by sending a pair $(gT,t)$ to the element $gtg^{-1}$ 
of $G$.  Clearly this is independent of the choice of 
representative of the coset $gT$.  This map has a number 
of very useful properties.  For instance it is $G$-equivariant 
for the obvious left $G$ action on $G/T\times T$ and 
the left $G$ action by conjugation on $G$.  More useful 
for us however is the following fact.  
\begin{proposition} 
\label{prop: Weyl map}
The map $p\colon G/T\times T\to G$ is a $|W|$ sheeted 
covering when restricted to the dense open subset $G_{\reg}$ 
of regular elements in $G$ (here $|W|$ denotes the order 
of the Weyl group $W$ of $G$).  
\end{proposition} 
Recall (see for instance \cite{Adams}) that an element $g$ of 
$G$ is said to be regular if the dimension of its centraliser 
is equal to $\dim T$.  
For $G = U(n)$ this means that $g$ has distinct eigenvalues or that it is conjugate to a 
diagonal matrix $\mathrm{diag}(z_1,\ldots, z_n)$ 
with all the $z_i$ distinct.  
It is clear that the set $G_{\reg}$ of 
regular elements in $G$ is a dense open subset of $G$.  

Over $G_{\reg}$ the map $p$ restricts to a map $p_{\reg}\colon 
G/T\times T_{\reg}\to G_{\reg}$ where $T_{\reg}$ denotes the 
dense open subset of $T$ consisting of regular elements.  
As mentioned in the statement of Proposition~\ref{prop: Weyl map} 
the map $p_{\reg}$ is a $|W|$ sheeted covering.  This means in particular 
the derivative $dp_{\reg}$ is surjective and so the pullback map 
$(p_{\reg})^*\colon \Omega(G_{\reg})\to \Omega(G/T\times T_{\reg})$ 
is injective on forms.  However since $G_{\reg}$ is dense 
in $G$ it follows that the pullback map 
$p^*\colon \Omega(G)\to \Omega(G/T\times T)$ must also be injective.  
We record the above discussion in the following Corollary to 
Proposition~\ref{prop: Weyl map}: 
\begin{corollary} 
\label{corr: f^* injective on forms}
The Weyl map $P \colon G/T \times T \to G$ induces by pullback an injective  map on 
differential forms
$$ 
p^*\colon \Omega(G)\to \Omega(G/T\times T).
$$  
\end{corollary}

For the special case when $G = U(n)$ and $T$ is the diagonal torus the map $p$ has a nice 
description. First of all $G/T$ can be identified 
with the flag manifold $F_n$ of $\CC^n$, so that we can think of a point 
in $G/T$ as an increasing sequence 
$$
V_1\subset V_2 \subset \cdots 
\subset V_n = \CC^n
$$ 
of subspaces of $\CC^n$ such that $\dim V_{i+1}/V_i = 1$,   
or alternatively as a family of 1-dimensional subspaces 
$W_1,W_2,\ldots, W_n$ 
of $\CC^n$ such 
that $W_i$ is orthogonal to $W_j$ if $i\neq j$.  Replacing the subspaces 
$W_i$ with the orthogonal projections $P_i$ onto them, we see that we 
can identify a point in $G/T$ with a family of orthogonal projections 
$P_1,P_2,\ldots, P_n$ such that $P_iP_j = 0$ if $i\neq j$ and $\sum_i P_i =1$.  

So we can think of a point in $G/T\times T$ as a pair $(P,\lambda)$ 
consisting of a family $P = (P_1,P_2,\ldots, P_n)$ of such orthogonal projections 
and a point $\lambda = (\lambda_1,\lambda_2,\ldots, \lambda_n)$ in $T$.  
It is useful to regard the $\lambda_i$ as the eigenvalues of a unitary 
matrix $g$ and the $P_i$ as the orthogonal projections onto the 
corresponding eigenspaces.  With this interpretation the map $p\colon G/T\times T\to
G$ has a simple description: it is the map which sends such a pair $(P,\lambda)$ 
to the unitary matrix 
$$ 
g = \sum^n_{i=1}\lambda_iP_i. 
$$ 
As we have already mentioned, 
the construction 
of the basic bundle gerbe in Section~\ref{sec:basic} gives 
in particular a bundle gerbe over the (diagonal) maximal 
torus $T = \TT^n$ of $G = U(n)$.  
Let us denote by $(L_T,Y_T)$ this basic bundle gerbe on 
$T$ so that $Y_T$ is the subset of $U(1)_0\times T$ consisting 
of pairs $(z,t)$ so that $z$ is not one of the diagonal entries of $t$.  
Note that the map $p\colon G/T\times T\to G$ induces a map 
\begin{gather*} 
p_Y\colon G/T\times Y_T\to Y \\ 
(gT,(z,t))\mapsto (z,gtg^{-1}).  
\end{gather*}
Let $Y_{\reg}$ be the subset of $Y$ consisting of pairs $(z,g)$ where 
$g\in G_{\reg}$.  Similarly let $Y_{T,\reg}$ denote the subset of $Y_T$ 
consisting of pairs $(z,t)$ where $t\in T_{\reg}$.  $Y_{\reg}$ and $Y_{T,\reg}$ 
are dense, open subsets of $Y$ and $Y_T$ respectively.  Also, the covering 
map $p_{\reg}\colon G/T\times T_{\reg}\to G_{\reg}$ pulls back along the 
projection $Y_{\reg}\to G_{\reg}$ to define a covering map $G/T\times Y_{T,\reg}\to Y_{\reg}$.  
This is just the restriction of the map $p_Y$ defined above.  Now 
the same argument used to prove Corollary~\ref{corr: f^* injective on forms} applies 
to prove the following Lemma.  
\begin{lemma} 
\label{lemm: f^* inj on forms} 
The map $p_Y\colon G/T\times Y_T\to Y$ induces by pullback an injective map on differential forms 
$$ 
(p_Y)^*\colon \Omega(Y)\to \Omega(G/T\times Y_T).  
$$ 
\end{lemma} 
One can prove in exactly the same way that the map $G/T\times Y_T^{[2]}\to Y^{[2]}$ 
induces an injective pullback map on differential forms.  Here the fibre products $Y_T^{[2]}$ and 
$Y^{[2]}$ are formed with respect to the projections $Y_T\to T$ and $Y\to G$ respectively.  

For later use, we will calculate the pullback $p^*(g^{-1}dg)$ of the 
operator valued 1-form $g^{-1}dg$ to $G/T\times Y_T$.  We get 
\begin{equation} 
\label{eq: pullback of MC form} 
p^*(g^{-1}dg) = \sum_{\lambda_i} \lambda_i^{-1}d\lambda_i P_i + 
\sum_{\lambda_i,\lambda_j}\lambda_i^{-1}\lambda_j P_i dP_j. 
\end{equation}

We finish this section with a remark about the eigenvalues of a unitary matrix.  
Clearly any unitary matrix $g$ can be written in the form $\sum_i \lambda_i P_i$ 
and it is even true that the eigenvalues $\lambda_i$ depend continuously on 
$g$, since they are the solutions of the equation $\det(g - \lambda 1) = 0$ and 
hence vary continuously with $g$.  
Consider the open subset $U_z$ of $G$ consisting of all unitary matrices 
$g$ such that $z$ is not an eigenvalue of $g$.  We could then define a 
partial order on the unit circle $U(1)\setminus \{z\}$ in the same manner as 
in Section~\ref{sec:basic} above.  We can then order the eigenvalues 
$\lambda_i(g)$ of the unitary matrices $g$ belonging to $U_z$ and we may 
wonder whether the $\lambda_i$ depend \emph{smoothly} on $g$.  
However this is not the case as the following example shows.  Take 
$G = SU(2)$ and consider the open set $U_i$ consisting of those 
$g\in SU(2)$ for which $i$ is not an eigenvalue.  Consider the 
intersection 
$$ 
U_i\cap T = \left\{\left(\begin{array}{cc} \a & 0 \\ 0 & \bar{\a}\end{array}\right)|\ \a\neq i\right\} 
$$ 
of $U_i$ with the standard diagonal torus $T$.  Let $g_{\a}$ denote the 
diagonal matrix above.  Then if we define $\lambda_1\colon U_1\to U(1)$ 
to be the first eigenvalue of $g_{\a}$ relative to the ordering of $U(1)\setminus \{i\}$ 
defined by rotating in a clockwise direction from $i$ then the value of 
$\lambda_1$ on $U_i\cap T$ is 
$$ 
\lambda_1(g_{\a}) = \begin{cases} 
\a & \text{if}\ x>0,y>0 \\ 
\bar{\a} & \text{if}\ x<0,y>0 \\ 
\a & \text{if}\ x<0,y<0\\ 
\bar{\a} & \text{if}\ x>0,y<0 
\end{cases} 
$$ 
with $\a = x+iy$.  This is continuous but not differentiable 
at $\a = \pm 1$.  In order to evaluate the various contour integrals we 
consider it is necessary for us to be able to write $g$ in the form 
$\sum_i \lambda_i P_i$.  Passing to the space $G/T\times T$ 
avoids this potential problem with the lack of smooth dependence on $g$  
of the eigenvalues $\lambda_i(g)$.

\section{The basic bundle gerbe as an equivariant bundle gerbe} 

Suppose that a compact Lie group 
$K$ acts smoothly on a manifold $M$, and that $(L,Y)$ is a bundle gerbe over $M$. 
\begin{definition} 
\label{def: equiv bundle gerbe} 
We will say that 
$(L,Y)$ is a \emph{$K$-equivariant bundle gerbe}\footnote{It is 
possible to consider a weaker notion of $K$-equivariant bundle 
gerbe --- see for instance \cite{Meinrenken}. The simpler notion 
that we describe here will be sufficient for our purposes.} if the 
following conditions are satisfied: 
\begin{enumerate} 
\item  there is an action of 
$K$ on $Y$ for which the surjective submersion $\pi\colon Y\to M$ 
is a $K$-equivariant map, 

\item the line bundle $L\to Y^{[2]}$ is a $K$-equivariant 
line bundle for the induced $K$-action on $Y^{[2]}$.  

\item the section $s$ of the line bundle $\delta(L)$ on $Y^{[3]}$ 
is $K$-equivariant. 
\end{enumerate} 
\end{definition} 
\noindent For more information 
on equivariant gerbes the reader is referred to \cite{BX,Meinrenken,Stienon}.

In the introduction we discussed the various realisations of the basic bundle 
gerbe that have appeared in the literature so far.  In the constructions 
\cite{BXZ, Bry1, Meinrenken} it is proven that the basic bundle gerbe is 
an equivariant bundle gerbe for the action of $G$ on itself by conjugation.  This is also true for our realisation of the basic 
bundle gerbe.  Let $G$ denote one of the unitary groups described in Section~\ref{sec:basic}.   
\begin{proposition} 
\label{prop: equiv} 
Let $G$ act on itself by conjugation.  Then the  basic bundle gerbe on $G$ constructed in Section~\ref{sec:basic} is an 
equivariant bundle gerbe in the sense of Definition~\ref{def: equiv bundle gerbe} 
above for this $G$ action. 
\end{proposition}

We first note that the conjugation action on $G$ lifts to an action on $Y$: if 
$(z,g)\in Y$ and $k\in G$ then $(z,kgk^{-1})\in Y$ since conjugation 
does not change the eigenvalues of a unitary operator.  This $G$-action naturally 
induces one on $Y^{[2]}$ and we need to check that the line bundle $L$ is 
equivariant for this $G$-action.  
Recall that $L\to Y^{[2]}_+$ is defined to be the top exterior power 
$$ 
L = \det(E) 
$$ 
where $E$ is the vector bundle on $Y^{[2]}_+$ defined by the projection valued 
map 
$$ 
P\colon Y^{[2]}\to B(H)
$$ 
$G$ acts naturally on $H$ and hence on $B(H)$; if $X$ is a bounded operator on $H$, 
then $k(X)$ is the bounded operator given by $k(X)(v) = kX(k^{-1}v)$ for $v$ a vector in $H$.  
Therefore, to prove that $L\to Y^{[2]}_+$ is equivariant, it is enough to prove that $E$ is 
an equivariant vector bundle, and to do that it is enough to prove that $P$ is an equivariant map.  

The value of $P$ at a point $(z_1,z_2,g)$ is the orthogonal projection onto the subspace 
$$ 
E_{(z_1,z_2,g)} = \bigoplus_{z_1 > \lambda > z_2} E_{(g,\lambda)} 
$$ 
If $v$ is an eigenvector of $g$ with eigenvalue $\lambda$ then clearly 
$k(v)$ is an eigenvector of $kgk^{-1}$ with eigenvalue 
$\lambda$.  Since $G$ acts as a group of unitary operators on $H$, it follows that 
$P(z_1,z_2,kgk^{-1}) = kP(z_1,z_2,g)k^{-1}$, i.e $P$ is equivariant.  
It is now clear how to extend the $G$-action on the line bundle $L$ over $Y^{[2]}_+$ 
to a $G$-action on $L$ over the entire space $Y^{[2]}$.  
It is also easy to see that the section $s$ of $\delta(L)$ defining the 
bundle gerbe product on $L$ is equivariant for these $G$-actions.  
This completes the proof.  

Notice in particular that Proposition~\ref{prop: equiv} implies that if 
$T$ is the diagonal torus inside $G = U(n)$ then $(L_T,Y_T)$ is a 
$T$-equivariant bundle gerbe for the trivial action of $T$ on itself 
by conjugation.  It is interesting to study the pullback of $(L,Y)$ 
along the map $p\colon G/T\times T$; in some sense the bundle gerbe 
$(L,Y)$ `abelianises'.  To understand what we mean by this observe that 
since $T$ is a subgroup of 
$G$ there is a canonical way to extend the $T$-equivariant bundle gerbe $(L_T,Y_T)$ 
on $T$ to a $G$-equivariant gerbe on $G/T\times T$.  We make 
$Y_T$ into a $G$-space by forming $G\times^T Y_T = G/T\times Y_T$, and we make 
$L_T$ into a $G$-equivariant line bundle over $G/T\times Y_T^{[2]}$ by 
forming $G\times^T L_T$.  The pair $(G\times^T L_T,G/T\times Y_T)$ is a 
$G$-equivariant bundle gerbe on $G/T\times T$.  We have the following proposition.  
\begin{proposition} 
There is a bundle gerbe morphism 
$$ 
(G\times^T L_T,G/T\times Y_T,G/T\times T)\to (L,Y,G) 
$$ 
covering the map $p\colon G/T\times T\to G$.  In particular the
bundle gerbe $p^*L$ on $G/T\times T$ has the same Dixmier-Douady class as 
$(G\times^T L_T,G/T\times Y_T)$.
\end{proposition} 
Here by a morphism of bundle gerbes we understand a morphism 
in the sense of \cite{Mur96Bundle-gerbes}.   
As we have already noted above, there 
is a canonical $G$-equivariant map 
$p_Y\colon G/T\times Y_T\to Y$ covering the Weyl map $p\colon G/T\times T\to G$,   
which sends a pair $(gT,(z,t))\in G/T\times Y_T$ to the pair $p_Y(gT,(z,t)) = (z,gtg^{-1})\in Y$. 

There is also a canonical $G$-equivariant map $G\times^T L_T\to L$ 
covering the induced map $p_Y^{[2]}\colon G/T\times Y_T^{[2]}\to Y^{[2]}$.  This is defined 
as follows.  Suppose that $(z_1,z_2,t)\in Y_{T,+}^{[2]}$, where $Y^{[2]}_{T,+}$ 
is defined in the analogous manner to $Y^{[2]}_+$.  Let $v\in E_{(z_1,z_2,t)}$ 
be an eigenvector of $t$ corresponding to some eigenvalue $\lambda \in (z_1,z_2)$.  
If $g\in G$ then $gv$ is an eigenvector of $gtg^{-1}$ corresponding to the eigenvalue 
$\lambda$.  In other words $gv\in E_{(z_1,z_2,gtg^{-1})}$.  We have therefore 
a linear map 
$$ 
E_{(z_1,z_2,t)} \to E_{(z_1,z_2,gtg^{-1})}. 
$$ 
On taking top exterior powers this gives a linear map 
$$ 
(L_T)_{(z_1,z_2,t)}\to L_{(z_1,z_2,gtg^{-1})}. 
$$ 
These linear maps are the restrictions to the fibres of  
a morphism of line bundles $\hat{p}\colon G\times^T L_T\to L$, covering 
$p_Y^{[2]}\colon G\times^T Y^{[2]}_{T,+}\to Y^{[2]}_+$.  It is clear that 
moreover this morphism is $G$-equivariant.  By duality we get a corresponding 
$G$-equivariant morphism of line bundles covering $p_Y^{[2]}\colon G\times^T Y^{[2]}_{T,-}\to Y^{[2]}_-$
Trivially a similar statement is true over $Y^{[2]}_{T,0}$.  Hence we 
have a $G$-equivariant morphism of $G$-equivariant line bundles 
$\hat{p}\colon G\times^T L_T\to L$ which covers $p_Y^{[2]}$.     
It is not difficult to show that $\hat{p}$ respects the bundle gerbe products 
on $G\times^T L_T$ and $L$.  The triple 
$$
(\hat{p},p_Y,p)\colon (G\times^T L_T,G/T\times Y_T,G/T\times T) \to (L,Y,G) 
$$ 
is a 
morphism of bundle gerbes in the sense of \cite{Mur96Bundle-gerbes}.  

\appendix

\section{Proof of Proposition~\ref{prop: expression for bg curvature}}  
\label{sec: first proof}
As in Section~\ref{sec:basic} 
above, if $(z_1,z_2,g)\in Y^{[2]}_+$ we let $P= P_{(z_1,z_2,g)}$ denote the 
orthogonal projection onto the eigenspaces $E_{(g,\lambda)}$ for 
$\lambda\in (z_1,z_2)$.  We have the contour integral expression~\eqref{eq: reg expression for P} for the restriction of $P$ 
to the open set $U_{(w_1,w_2,h)} = U_{w_1}\times U_{w_2}\times U_h$ associated to a point 
$(w_1,w_2,h)\in Y^{[2]}_+$: 
$$ 
P = \frac{1}{2\pi i}\oint_{C}(\xi - g)^{-1}d\xi 
$$ 
Here $C$ is the contour in $U_{(w_1,w_2,h)}$ described earlier. 
This formula clearly shows that $P$ is differentiable, since the 
resolvent $(\xi - g)^{-1}$ depends smoothly on $g$.  It is not hard to show that 
the derivative $dP$ of $P$ at a point $(z_1,z_2,g)$ in the open set $U_{(w_1,w_2,h)}$ is given by the expression
$$ 
dP = \frac{1}{2\pi i}\oint_{C}(\xi - g)^{-1}dg(\xi - g)^{-1}d\xi 
$$ 
where as above, $dg$ denotes the derivative at $(z_1,z_2,g)$ of the projection of $Y^{[2]}$ onto $G$.  
Note that the integrand is an analytic function of $\xi$ and hence we can deform 
the contour $C$ without changing the value of the integral.  In particular 
we may write 
$$ 
dP = \frac{1}{2\pi i}\oint_{C_{(z_1,z_2,g)}}(\xi - g)^{-1}dg(\xi - g)^{-1}d\xi 
$$ 
where for any $(z_1,z_2,g)\in Y^{[2]}_+$, the contour $C_{(z_1,z_2,g)}$ is chosen so that it surrounds the part of 
the spectrum of $g$ lying between $z_1$ and $z_2$.   
We now 
compute the expression for the curvature 
$$
F_{\det(\nabla)} = \tr(PdPdP)
$$ 
in terms of these contour integral formulas.  First of all choose contours $C_{(z_1,z_2,g)}$, $C'_{(z_1,z_2,g)}$ 
and $C''_{(z_1,z_2,g)}$ 
surrounding the part of the spectrum of $g$ lying in the set $(z_1,z_2)$.  We can suppose that the contours are 
nested in the sense that $C_{(z_1,z_2,g)}$ is contained inside $C'_{(z_1,z_2,g)}$ which is contained inside 
$C''_{(z_1,z_2,g)}$.    
Then $F_{\det(\nabla)}$ can be 
written as an iterated contour integral 
$$
\left(\frac{1}{2\pi i}\right)^3\oint_{C_{(z_1,z_2,g)}}\oint_{C'_{(z_1,z_2,g)}}\oint_{C''_{(z_1,z_2,g)}} 
\tr\left( R(\xi,\eta,\zeta,g)\right) d\xi d\eta d\zeta
$$ 
where $R(\xi,\eta,\zeta,g)$ is the product of resolvents 
$$ 
R(\xi,\eta,\zeta,g) = 
(\xi-g)^{-1}(\eta-g)^{-1}dg(\eta-g)^{-1} 
(\zeta - g)^{-1}dg 
(\zeta - g)^{-1}
$$
Using the resolvent formula we can rewrite the products appearing in the 
integrand as differences, for example  
$$ 
(\xi-g)^{-1}(\eta-g)^{-1} = (\eta -\xi)\left[ (\xi - g)^{-1} - (\eta - g)^{-1}\right]. 
$$ 
Using this fact 
and the cyclic property of the trace we can rewrite the integrand as 
$(\zeta - \eta)^{-1}(\zeta - \xi)^{-1}(\eta - \xi)^{-1}$ times the 
two-form 
$$
\tr\left( 
\left[ (\xi - g)^{-1} - (\eta - g)^{-1}\right]dg\left[(\eta - g)^{-1} - (\zeta - g)^{-1}\right]dg\right) 
$$
Expanding out the product inside the trace leaves us with 
four contour integrals to compute.  Of these, the only non-vanishing contour integral 
gives    
$$
\left(\frac{1}{2\pi i}\right)^3\oint_{C_{(z_1,z_2,g)}}\oint_{C'_{(z_1,z_2,g)}}\oint_{C''_{(z_1,z_2,g)}} 
\frac{\tr\left((\xi - g)^{-1}dg(\zeta - g)^{-1}dg\right)}{(\zeta - \eta)(\zeta - \xi)(\eta - \xi)}\,  d\xi d\eta d\zeta 
$$
which we can evaluate to 
\begin{equation} 
\label{curv expression} 
\left(\frac{1}{2\pi i}\right)^2\oint_{C_{(z_1,z_2,g)}}\oint_{C''_{(z_1,z_2,g)}}(\zeta - \xi)^{-2}\tr\left((\xi - g)^{-1}dg 
(\zeta - g)^{-1}dg\right)d\xi d\zeta. 
\end{equation}
The function $(\zeta - \xi)^{-2}$ is a holomorphic function of $\xi$, since 
$\zeta \in C''_{(z_1,z_2,g)}$ which lies outside $C_{(z_1,z_2,g)}$.  Therefore we can simplify this 
expression further using the holomorphic functional calculus.  
Since $(\zeta - \xi)^2(\zeta - \xi)^{-2} = 1$ we have (using the property 
$f(T)f'(T) = (ff')(T)$ of the functional calculus) 
$$ 
\frac{1}{2\pi i}\oint_{C_{(z_1,z_2,g)}}(\zeta - \xi)^{-2}(\xi - g)^{-1}d\xi = (\zeta - g)^{-2}P_{z_1z_2}. 
$$ 
Therefore, after renaming variables, we can rewrite~\eqref{curv expression} above as the single contour 
integral 
$$ 
\frac{1}{2\pi i}\oint_{C_{(z_1,z_2,g)}} \tr\left( (\xi - g)^{-2}P_{z_1z_2}dg(\xi - g)^{-1}dg\right) d\xi. 
$$
We want to show that we have the equality of 
Proposition~\ref{prop: expression for bg curvature}:
\begin{multline}
\label{eq: equality of two-forms} 
\frac{1}{2\pi i}\oint_{C_{(z_1,z_2,g)}} \tr\left((\xi - g)^{-1}dg(\xi - g)^{-2}P_{z_1z_2}dg\right) d\xi \\ 
= \frac{1}{4\pi i}\oint_{C_{(z_1,z_2,g)}} \tr\left((\xi - g)^{-1}dg(\xi - g)^{-2}dg\right) d\xi
\end{multline} 
Up to now everything that we have said works equally well 
for $G$ finite dimensional or for $G$ one of the infinite dimensional 
groups $U_p(H)$.  Let us now suppose that the group 
$G$ is the finite dimensional unitary group $U(H)$, 
for $H$ a finite dimensional complex Hilbert space.  
Thus $G$ is isomorphic to $U(n)$ where $n$ is the 
dimension of $H$.  
 This is the first instance where we will take advantage of the 
special properties of the map $p\colon G/T\times T\to G$.  From 
the remark following Lemma~\ref{lemm: f^* inj on forms} we have that 
the pullback map $\Omega(Y^{[2]})\to \Omega(G/T\times Y_T)$ 
induced by $G/T\times Y_T^{[2]}\to Y^{[2]}$ is injective.  
If we think of a point of $G/T\times Y_T^{[2]}$ as a family $(P_i,(z_1,z_2,\lambda_i))$ 
then the map $G/T\times Y_T^{[2]}\to Y^{[2]}$ can be written 
as 
$$ 
(P_i,(z_1,z_2,\lambda_i))\mapsto (z_1,z_2,\sum \lambda_iP_i) 
$$
Using this description of the map $G/T\times Y_T^{[2]}\to Y^{[2]}$ 
in terms of the projections $P_i$ and 
eigenvalues $\lambda_i$ we see that we 
can write the pullback by $p$ of the left hand side of~\eqref{eq: equality of two-forms} as 
$$ 
\frac{1}{2\pi i}\oint_{C_{(z_1,z_2,g)}} \sum_{i,j} (\xi - \lambda_i)^{-1} 
(\xi -\lambda_j)^{-2}\tr\left(P_i p^*(dg)P_j P_{z_1z_2}p^*(dg)\right) d\xi
$$ 
To simplify this expression we will make use of the following 
easily proven facts: 
\begin{align*} 
& \Res_{\xi =\lambda_j}(\xi - \lambda_i)^{-1}(\xi - \lambda_j)^{-2} = - (\lambda_i - \lambda_j)^{-2} \\  
& \Res_{\xi =\lambda_i}(\xi - \lambda_i)^{-1}(\xi - \lambda_j)^{-2} = (\lambda_i - \lambda_j)^{-2}\\ 
& P_jP_{z_1z_2} = \begin{cases} 
P_j & \text{if}\ j\in (z_1,z_2) \\ 
0 & \text{if}\ j\notin (z_1,z_2) 
\end{cases} 
\end{align*}
By splitting the sum over the eigenvalues $\lambda_i$ into the sum over the eigenvalues $\lambda_i\in (z_1,z_2)$ 
and the eigenvalues $\lambda_i\notin (z_1,z_2)$ and using the above facts, we derive the 
following expression for the pullback by $p$ of the left hand side of~\eqref{eq: equality of two-forms}: 
$$ 
- \sum_{\stackrel{\scriptstyle{\lambda_i \notin (z_1,z_2)}}{ 
\lambda_j \in (z_1,z_2)}}(\lambda_i - \lambda_j)^{-2}\tr(P_ip^*(dg)P_jp^*(dg)) 
$$
On the other hand, it is easy to see that the pullback by $p$ of the right hand side of~\eqref{eq: equality of two-forms} can be written as 
$$ 
\frac{1}{4\pi i}\oint_{C_{(z_1,z_2,g)}} \sum_{i,j}(\xi - \lambda_i)^{-1}(\xi - \lambda_j)^{-2} 
\tr(P_i p^*(dg) P_j p^*(dg)) d\xi
$$ 
Again, we can split the sum over eigenvalues $\lambda_i$ and $\lambda_j$ into 
sums over eigenvalues belonging to the sets $(z_1,z_2)$ or their complements to get 
\begin{align*}
& \frac{1}{4\pi i}\oint_{C_{(z_1,z_2,g)}} \sum_{ 
\stackrel{\scriptstyle{\lambda_i \in (z_1,z_2)}}{\lambda_j \in (z_1,z_2)} } (\xi - \lambda_i)^{-1}(\xi - \lambda_j)^{-2} 
\tr(P_i p^*(dg) P_j p^*(dg)) d\xi \\ 
& + \frac{1}{4\pi i}\oint_{C_{(z_1,z_2,g)}} \sum_{ 
\stackrel{\scriptstyle{\lambda_i \notin (z_1,z_2)}}{\lambda_j \in (z_1,z_2)} } (\xi - \lambda_i)^{-1}(\xi - \lambda_j)^{-2} 
\tr(P_i p^*(dg) P_j p^*(dg)) d\xi \\ 
& + \frac{1}{4\pi i}\oint_{C_{(z_1,z_2,g)}} \sum_{ 
\stackrel{\scriptstyle{\lambda_i \in (z_1,z_2)}}{\lambda_j \notin (z_1,z_2)} } (\xi - \lambda_i)^{-1}(\xi - \lambda_j)^{-2} 
\tr(P_i p^*(dg) P_j p^*(dg)) d\xi 
\end{align*} 
By the residue theorem these contour integrals become 
\begin{align*} 
- \frac{1}{2}  \sum_{ 
\stackrel{\scriptstyle{\lambda_i \notin (z_1,z_2)}}{\lambda_j \in (z_1,z_2)} } 
(\lambda_i - \lambda_j)^{-2}\tr\left( P_i p^*(dg) P_j p^*(dg)\right) \\ 
+ \frac{1}{2}  \sum_{ 
\stackrel{\scriptstyle{\lambda_i \in (z_1,z_2)}}{\lambda_j \notin (z_1,z_2)} } 
(\lambda_i - \lambda_j)^{-2}\tr\left( P_i p^*(dg) P_j p^*(dg)\right) \\ 
= - \sum_{ 
\stackrel{\scriptstyle{\lambda_i \notin (z_1,z_2)}}{\lambda_j \in (z_1,z_2)} } 
(\lambda_i - \lambda_j)^{-2} \tr\left( P_ip^*(dg) P_j p^*(dg)\right) 
\end{align*} 
and so we see the two expressions are equal.  It follows by the injectivity of the map 
$p^*\colon \Omega^2(Y^{[2]})\to \Omega^2(G/T\times Y_T^{[2]})$ that the two forms on $Y^{[2]}$ are equal.  

Now let us suppose that $G$ is one of the infinite dimensional Banach Lie 
groups $U_p(H)$ for $1\leq p\leq 2$.  We have a pair of two-forms 
defined on $Y^{[2]}_+$ 
and we want to prove that we have the equality of~\eqref{eq: equality of two-forms}.  Since neither of these 
two-forms have components in the $U(1)$ directions, it is sufficient to prove 
that the two-forms are equal for fixed $z_1$ and $z_2$.  Thus we can regard them 
as forms on the open subset $U_{z_1z_2}$ of $G$ consisting of those $g\in G$ for which neither 
$z_1$ nor $z_2$ is an eigenvalue.  To prove that they are equal, we will use a slight modification 
of the argument used to prove Theorem~\ref{thm: Quillens thm} and conclude that if $\omega$ is 
a form defined on an open subset $U$ of $G$ such that $\omega$ vanishes on restriction 
to $U\cap U(H_n)$ for all $n$, then $\omega$ is zero. The equality of~\eqref{eq: equality of two-forms} 
follows easily from this conclusion.  

Recall that Quillen proves Theorem~\ref{thm: Quillens thm} using the `tame approximation' 
theorem of Palais \cite{Palais}, which states that for any $A\in \cL_p$, 
$P_nAP_n\to A$ in $\cL_p$.  A $p$-form $\omega$ on $GL_p(H)$ can 
be regarded as a smooth map $\omega(A_1,\ldots, A_p)(g)$ which is 
alternating and multilinear in the $A_i\in \cL_p$.  If $\omega$ vanishes 
on restriction to $GL(H_n)$ for every $n$ then it follows from the tame approximation theorem 
that $\omega$ vanishes on $GL_p(H)$.  More generally if $U$ is an 
open subset of $GL_p(H)$ and $\omega$ is a $p$-form on $U$ which 
vanishes on restriction to any $GL(H_n)\cap U$ for every $n$ then $\omega = 0$.  
To prove the theorem for the unitary groups $U_p(H)$ Quillen considers 
the phase retraction map 
\begin{align*} 
& GL_p(H)\to U_p(H) \\  
& g\mapsto g(gg^*)^{-1/2}. 
\end{align*} 
It is easy to see that this maps the subgroups $GL(H_n)$ of $GL_p(H)$ 
to the corresponding subgroups $U(H_n)$ of $U_p(H)$.  The phase 
retraction has the property that it maps forms on $U_p(H)$ injectively 
to forms on $GL_p(H)$.  The result for $U_p(H)$ now follows from the 
corresponding result for $GL_p(H)$.  Again let us note that if $\omega$ 
is a $p$-form defined on an open subset $U\subset U_p(H)$ such that 
$\omega$ vanishes on restriction to any $U(H_n)\cap U$, then 
$\omega$ vanishes on $U$.  To see this, note that if $V$ is the inverse 
image of $U$ under the phase retraction map, then the inverse image 
of $U(H_n)\cap U$ is mapped to $GL(H_n)\cap V$.  If $\omega'$ 
denotes the pullback of $\omega$ to $V$, then $\omega'$ vanishes 
on restriction to any $GL(H_n)\cap V$.  Therefore $\omega'$, and 
hence $\omega$, is zero.    

%

\section{Proof of Theorem~\ref{thm: main thm} part ({\rm b}).}

The map $p_Y\colon G/T\times Y_T\to Y$ from Lemma~\ref{lemm: f^* inj on forms} 
fits into a commutative diagram 
$$ 
\xymatrix{ 
G/T\times Y_T \ar[d]_-{p_Y} \ar[r] & G/T\times T \ar[d]^-p \\ 
Y \ar[r]^-\pi & G } 
$$ 
where we recall that we write $\pi$ for the projection map $Y\to G$.  
We want to show that $df = \pi^*(2\pi i\nu)$ where $\nu$ is the basic three-form~\eqref{eq: basic 3-form on SU(n)}.  It is sufficient to show this equality on the dense open 
subset $Y_{\reg}\subset Y$.  $(p_Y)^*$ is then injective on $\Omega(Y_{\reg})$ since 
$p_Y\colon G/T\times Y_{T,\reg}\to Y_{\reg}$ is a covering map.  
Note that by 
commutativity of the diagram, $(p_Y)^*\pi^*(2\pi i\nu)$ is equal to the pullback of $p^*(2\pi i\nu)$ along 
the map $G/T\times Y_T\to G/T\times T$.  We begin with the left hand side of the equation 
$(p_Y)^*df = (p_Y)^*\pi^*(2\pi i\nu)$.   
We  first calculate the 
pullback $(p_Y)^*f$ of the restriction of $f$ to $Y_{\reg}$.  Note that a point in $G/T\times Y_{T,\reg}$ is 
of the form $((P,\lambda),z)$ where $(P,\lambda) = (P_i,\lambda_i)$ is a family 
of projections $P_i$ and eigenvalues $\lambda_i$ and where $z$ is not equal to any of the $\lambda_i$.  
Since $(P,\lambda)\in G/T\times T_{\reg}$, all of the $\lambda_i$ are {\em distinct}.  Thus if $i\neq j$ then 
$P_iP_j = P_jP_i = 0$.  The pullback $(p_Y)^*f$ is given by the contour integral 
\begin{multline*} 
\frac{1}{8\pi^2}\oint_{C_{(z,g)}} \sum \log_z \xi (\xi - \lambda_i)^{-1} 
(\xi - \lambda_k)^{-2}\tr(P_i(d\lambda_j P_j  + \\ 
+ \lambda_jdP_j)P_k(d\lambda_l P_l + 
\lambda_l dP_l))d\xi 
\end{multline*} 
We simplify the term within the trace  summed  over $j$ and $l$.  Expanding it out and reindexing gives us 
\begin{align*}
\sum_{j,l} & \tr(d\lambda_j \lambda_l \delta_{ij} \delta_{jk}P_jdP_l - \lambda_j d\lambda_i \delta_{kl}\delta_{li}  P_idP_j 
+ \lambda_j\lambda_lP_idP_jP_kdP_l) \\
&= \sum_{l}  \tr(d\lambda_i \lambda_l \delta_{ik} P_i dP_l) - \sum_{j} \tr(\lambda_j d\lambda_i \delta_{ik}  P_idP_j )
+ \sum_{j,l} \tr(\lambda_j\lambda_lP_idP_jP_kdP_l) \\
&=\sum_{j} \tr(\lambda_j\lambda_lP_idP_jP_kdP_l) \\
\end{align*}
Inserting this expression back into the contour integral gives us
\begin{equation}
\label{eq: expression 1}
\frac{1}{8\pi^2}\oint_{C_{(z,g)}} \sum \log_z\xi   (\xi - \lambda_i)^{-1}(\xi - \lambda_k)^{-2}
\lambda_j\lambda_l\tr(P_idP_jP_kdP_l)d\xi. 
\end{equation}
To evaluate this contour integral we use Cauchy's residue theorem.  For each pair of indices $i$ and 
$k$ we need to calculate the residues of the function $\log_z\xi (\xi - \lambda_i)^{-1}(\xi - \lambda_k)^{-2}$ 
at the poles $\xi = \lambda_i$ and $\xi = \lambda_k$.  When $i\neq k$ we have  
$$ 
\Res_{\xi = \lambda_i}(\xi -\lambda_i)^{-1}(\xi - \lambda_k)^{-2} \log_z\xi = 
\log_z\lambda_i (\lambda_i - \lambda_k)^{-2} 
$$ 
and 
$$ 
\Res_{\xi = \lambda_k} (\xi - \lambda_i)^{-1}(\xi - \lambda_k)^{-2} \log_z\xi = 
(\lambda_k - \lambda_i)^{-1}\lambda_k^{-1} - \log_z \lambda_k (\lambda_k - 
\lambda_i)^{-2}. 
$$ 
We also need to consider the residue at the pole $\xi = \lambda_i = \lambda_k$ in the 
case where $i=k$.  In this case however, the two-form $\tr(P_idP_jP_idP_k)$ vanishes.  
This is clear when $j=i$ as $P_idP_iP_i = 0$.  When $j\neq i$ we can use the identity 
$dP_iP_j + P_idP_j = 0$ to write $P_idP_jP_i = - dP_iP_jP_i = 0$.  Therefore, only the residues 
at the poles $\xi = \lambda_i$ and $\xi = \lambda_k$ for $i$ and $k$ distinct give contributions.      
Therefore~\eqref{eq: expression 1} becomes 
\begin{multline} 
\label{eq: expression 2} 
\frac{i}{4\pi } \sum_{i\neq k} \left( 
\log_z \lambda_i (\lambda_i - \lambda_k)^{-2} - 
\log_z \lambda_k (\lambda_i - \lambda_k)^{-2} \right. \\
\left. + (\lambda_k - \lambda_i)^{-1}\lambda_i^{-1}\right) 
\lambda_j \lambda_l \tr(P_idP_j P_k dP_l) 
\end{multline} 
To simplify this expression we consider, for fixed 
$\lambda_i$ and $\lambda_k$, the sum 
\begin{equation} 
\label{eq: expression 3} 
\sum_{j,l}\lambda_j\lambda_l \tr(P_i dP_j P_k dP_l) 
\end{equation}
To evaluate this sum we shall make use of the following identity:
$$ 
\sum_j \lambda_j P_idP_j P_k = \lambda_i P_idP_iP_k + \lambda_k P_i dP_k P_k. 
$$ 
To see this note we write the sum as 
$$
\sum_j \lambda_j P_idP_jP_k = \lambda_i P_idP_iP_k + \lambda_k P_idP_kP_k + \sum_{i\neq j, j\neq k} \lambda_j P_idP_jP_k 
$$
If the indices $i$, $j$ and $k$ are such that $i\neq j$ and $j\neq k$, then we 
can write 
$$ 
P_idP_jP_k = -dP_iP_jP_k = 0, 
$$ 
and hence the identity.   
Applying this identity we find that~\eqref{eq: expression 3}
becomes a sum of two terms: 
$$ 
\sum_{l} \lambda_i\lambda_l \tr(P_idP_iP_k dP_l) 
+ \sum_{l} \lambda_k\lambda_l\tr(P_idP_kP_kdP_l) 
$$ 
We can apply the identity again to each of the terms in this new sum, after first writing 
$\tr(P_idP_iP_k dP_l) = \tr(P_idP_iP_k dP_lP_i)$ and 
$\tr(P_idP_kP_kdP_l) = \tr(P_idP_kP_kdP_lP_i)$.  We obtain 
\begin{multline*}
\lambda_i\lambda_k \tr(P_idP_iP_kdP_k) + \lambda_i^2\tr(P_idP_iP_kdP_i) \\
+ \lambda_k^2\tr(P_idP_kP_kdP_k) + \lambda_i\lambda_k \tr(P_idP_kP_kdP_i). 
\end{multline*}
Writing $dP_k = dP_k P_k + P_kdP_k$, $dP_i = dP_i P_i + P_idP_i$ 
and using the identities $P_iP_k = \delta_{ik}P_i$ and $dP_kP_i = \delta_{ik}dP_k - P_kdP_i$, we see that we can write this as 
\begin{multline*} 
(1- \delta_{ik})\lambda_i\lambda_k\tr(P_idP_idP_k) - 
(1-\delta_{ik})\lambda_i^2\tr(P_kdP_idP_i)  \\
+ (1-\delta_{ik})\lambda_k^2\tr(P_idP_kdP_k) 
- (1-\delta_{ik})\lambda_i\lambda_k\tr(P_idP_kdP_k). 
\end{multline*}
We can further simplify this, using the above properties of the projections 
$P_i$ and $P_k$, to 
$$ 
(1-\delta_{ik})(\lambda_i - \lambda_k)^2\tr(P_idP_kdP_k) = (\lambda_i - \lambda_k)^2\tr(P_idP_kdP_k). 
$$ 
Substituting this into~\eqref{eq: expression 2} gives the following 
expression for the pullback $(p_Y)^*f$:
\begin{equation}
\label{eq: expression 4}
\frac{i}{4\pi }\sum_{i\neq k}\left( 
\log_z \lambda_i - \log_z\lambda_k +(\lambda_k - \lambda_i)\lambda_k^{-1}\right) 
\tr(P_idP_kdP_k). 
\end{equation} 
On taking the exterior derivative we end up with the 
following expression for $d(p_Y)^*f$: 
\begin{multline} 
\label{eq: expression 5}
\frac{i}{4\pi }\sum_{i\neq k} 
\left( \lambda_i^{-1}d\lambda_i - \lambda_k^{-1}
d\lambda_k - d\lambda_i \lambda_k^{-1} 
+ \lambda_i \lambda_k^{-1}d\lambda_k \lambda_k^{-1} \right) 
\tr(P_idP_kdP_k) \\ 
+ \frac{i}{4\pi}\sum_{i\neq k} \left(\log_z \lambda_i - \log_z \lambda_k + 1 - \lambda_i\lambda_k^{-1}\right)\tr(dP_idP_kdP_k).  
\end{multline} 
We can simplify the second term in this expression further: since $\tr(dP_idP_idP_i) = 0$ we see that 
we can write $\sum_{i\neq k}\tr(dP_idP_kdP_k) = \sum_{i,k}\tr(dP_idP_kdP_k)$ 
and this last sum is zero since $\sum dP_i = 0$.  We can also reindex and write 
\begin{multline*} 
\sum_{i\neq k}\log_z \lambda_i\tr(dP_idP_kdP_k) - 
\sum_{i\neq k}\log_z\lambda_k\tr(dP_idP_kdP_k) \\ 
= \sum_{i\neq k}\log_z\lambda_i\tr(dP_idP_kdP_k) - 
\sum_{i\neq k}\log_z\lambda_i\tr(dP_kdP_idP_i) 
\end{multline*}
If $i\neq k$ then it is easy to see that $\tr(dP_idP_kdP_k) = - \tr(dP_kdP_idP_i)$.  
Therefore we can write the above expression as 
$$ 
-2\sum_{i\neq k} \log_z\lambda_i \tr(dP_kdP_idP_i). 
$$ 
Again, since $\tr(dP_idP_idP_i) = 0$ we can write this as 
$$ 
-2\sum_{i,k}\log_z\lambda_i\tr(dP_kdP_idP_i)
$$ 
which equals zero since $\sum dP_k = 0$.  Therefore we end up with the following expression 
for $(p_Y)^*df$: 
\begin{multline} 
\label{eq: expression for p*df} 
\frac{i}{4\pi }\sum_{i\neq k} 
\left( \lambda_i^{-1}d\lambda_i - \lambda_k^{-1}
d\lambda_k - d\lambda_i \lambda_k^{-1} 
+ \lambda_i \lambda_k^{-1}d\lambda_k \lambda_k^{-1} \right) 
\tr(P_idP_kdP_k) \\ 
- \frac{i}{4\pi}\sum_{i\neq k} \lambda_i\lambda_k^{-1}\tr(dP_idP_kdP_k). 
\end{multline}
We would now like to compare this expression with 
the pullback three-form $p^*(2\pi i\nu)$ where we remind the reader again that $\nu$ is the three-form~\eqref{eq: basic 3-form on SU(n)}.  To compute this 
recall from equation~\eqref{eq: pullback of MC form} that we had 
$$ 
p^*(g^{-1}dg) =  \lambda_i^{-1}d\lambda_i P_i + \lambda_i^{-1}\lambda_jP_idP_j. 
$$ 
After a little calculation one finds that the pullback three-form $2\pi i p^*\nu$ is given by 
the following sum of two terms  
\begin{equation}
\label{eq: pullback of nu}
-\frac{i}{4\pi}\tr\left(\lambda_i^{-1}d\lambda_i\lambda_i^{-1}
\lambda_j\lambda_k^{-1}\lambda_lP_idP_jP_kdP_l\right) 
- \frac{i}{12\pi}\frac{\lambda_j\lambda_l\lambda_n}{\lambda_i\lambda_k\lambda_m}\tr\left(P_idP_j 
P_kdP_lP_mdP_n\right)
\end{equation} 
where in each term there is understood to be a sum over all appropriate indices.
We concentrate on the first term in~\eqref{eq: pullback of nu} to begin with (for clarity 
we omit the factor of $-i/4\pi$).  
Making use of the fact that $dP_i P_j + P_idP_j = \delta_{ij}dP_i$ (so that $P_idP_j = dP_i(\delta_{ij} - P_j)$) 
we see that we can write it as 
\begin{align*} 
&\sum \lambda_i^{-1}d\lambda_i\lambda_i^{-1}\lambda_j\lambda_k^{-1}\lambda_l
\tr(P_idP_jP_kdP_l) \\ 
=& \sum \lambda_i^{-1}d\lambda_i\lambda_i^{-1}\lambda_j\lambda_k^{-1}\lambda_l 
\tr(\delta_{ij}dP_i - dP_iP_j)P_kdP_l) \\ 
= &\sum \left(\left[\lambda_i^{-1}d\lambda_i\lambda_k^{-1}\lambda_l 
- \lambda_i^{-1}d\lambda_i\lambda_i^{-1}\lambda_l\right]\tr(dP_iP_kdP_l)\right) \\ 
= & \sum \left(\left[\lambda_i^{-1}d\lambda_i\lambda_k^{-1}\lambda_l 
- \lambda_i^{-1}d\lambda_i\lambda_i^{-1}\lambda_l\right]\tr\left\{dP_i(\delta_{kl}dP_k - dP_kP_l)
\right\}\right) \\ 
= & \sum \left(\lambda_i^{-1}d\lambda_i\tr(dP_idP_k) - \lambda_i^{-1}d\lambda_i\lambda_k^{-1}
\lambda_l\tr(dP_idP_kP_l) \right. \\ 
& \left. - \lambda_i^{-1}d\lambda_i\lambda_i^{-1}\lambda_k\tr(dP_idP_k) 
+ \lambda_i^{-1}d\lambda_i\lambda_i^{-1}\lambda_l\tr(dP_idP_kP_l) \right) 
\end{align*}
We can write the terms $\lambda_i^{-1}d\lambda_i\lambda_k^{-1}\lambda_l\tr(dP_idP_kP_l)$ and 
$\lambda_i^{-1}d\lambda_i\lambda_i^{-1}\lambda_l\tr(dP_idP_kP_l)$ as 
\begin{multline*} 
d\lambda_i\lambda_k^{-1}\tr(dP_idP_kP_i) + \lambda_i^{-1}d\lambda_i\tr(dP_idP_kP_k) 
+ \sum_{i\neq l, k\neq l}\lambda_i^{-1}d\lambda_i\lambda_k^{-1}\lambda_l\tr(dP_idP_kP_l) 
\end{multline*} 
and 
\begin{multline*} 
\lambda_i^{-1}d\lambda_i\tr(dP_idP_kP_i) + \lambda_i^{-1}d\lambda_i\lambda_i^{-1}
\lambda_k\tr(dP_idP_kP_k)  \\ 
+ \sum_{i\neq l,k\neq l}\lambda_i^{-1}d\lambda_i\lambda_i^{-1}
\lambda_l\tr(dP_idP_kP_l). 
\end{multline*} 
respectively.  Note that if $i\neq l$ and $k\neq l$ then $\tr(dP_idP_kP_l) = -\delta_{ik}\tr(P_idP_ldP_l)$.  
This is because 
$$ 
\tr(dP_idP_kP_l) = 
\tr(P_ldP_idP_kP_l) = \tr(dP_lP_iP_kdP_l) 
$$ 
It follows that 
$$ 
\sum_{i\neq l,k\neq l}\left(\lambda_i^{-1}d\lambda_i\lambda_i^{-1}\lambda_l - 
\lambda_i^{-1}d\lambda_i\lambda_k^{-1}\lambda_l\right)\tr(dP_idP_kP_l) = 0. 
$$ 
Therefore we end up with the following expression: 
\begin{multline}
\label{eq: expression 6} 
\sum\left(\lambda_i^{-1}d\lambda_i\tr(dP_idP_k) - \lambda_i^{-1}d\lambda_i\tr(dP_idP_kP_k) 
- d\lambda_i\lambda_k^{-1}\tr(dP_idP_kP_i) \right. \\ 
\left. - \lambda_i^{-1}d\lambda_i\lambda_i^{-1}
\lambda_k\tr(dP_idP_k) + \lambda_i^{-1}d\lambda_i\tr(dP_idP_kP_i) \right. \\ 
\left. + \lambda_i^{-1}d\lambda_i\lambda_i^{-1}\lambda_k\tr(dP_idP_kP_k)\right). 
\end{multline} 
Using $dP_k = P_kdP_k + dP_kP_k$ we can simplify the terms 
$$ 
\sum \left( \lambda_i^{-1}d\lambda_i \tr(dP_idP_k) - 
\lambda_i^{-1}d\lambda_i \tr(dP_idP_kP_k) \right)
$$ 
appearing in~\eqref{eq: expression 6} to 
$$ 
-\sum \lambda_i^{-1}d\lambda_i \tr(P_kdP_k dP_i). 
$$ 
Similarly we can simplify the terms 
$$ 
\sum \left( \lambda_i^{-1}d\lambda_i\lambda_i^{-1}\lambda_k\tr(dP_idP_kP_k) 
- \lambda_i^{-1}d\lambda_i \lambda_i^{-1}\lambda_k \tr(dP_idP_k) \right)
$$ 
appearing in~\eqref{eq: expression 6} to 
$$ 
\sum \lambda_i^{-1}d\lambda_i \lambda_i^{-1}\lambda_k \tr(P_kdP_kdP_i). 
$$ 
Therefore~\eqref{eq: expression 6} becomes 
\begin{multline} 
\label{eq: expression 6.5} 
\sum\left( - \lambda_i^{-1}d\lambda_i\tr(P_kdP_kdP_i) - d\lambda_i\lambda_k^{-1}\tr(dP_idP_kP_i) \right. \\  
\left. + \lambda_i^{-1}d\lambda_i\lambda_i^{-1}\lambda_k\tr(P_kdP_kdP_i) + \lambda_i^{-1}d\lambda_i\tr(dP_idP_kP_i)\right) 
\end{multline} 
Replacing $dP_k$ with $dP_kP_k + P_kdP_k$ we can write 
$$
d\lambda_i\lambda_k^{-1}\tr(dP_idP_kP_i) = 
d\lambda_i\lambda_k^{-1}\tr(dP_iP_kdP_kP_i)
$$ 
which is then easily seen to be equal to 
$-d\lambda_i\lambda_k^{-1}\tr(P_idP_kdP_k)$.  Similarly write 
$$
\lambda_i^{-1}d\lambda_i\lambda_i^{-1}\lambda_k\tr(P_kdP_kdP_i) =  
\lambda_i^{-1}d\lambda_i\lambda_i^{-1}\lambda_k\tr(P_kdP_kP_idP_i)
$$ 
which is in turn equal 
to $-\lambda_i^{-1}d\lambda_i\lambda_i^{-1}\lambda_k\tr(P_kdP_idP_i)$.  Then after reindexing some of the sums 
in~\eqref{eq: expression 6.5} and restoring the factor $-i/4\pi$ we get 
\begin{equation}
\label{eq: first part of p*nu} 
\frac{i}{4\pi}\sum_{i\neq k}\left( 
\lambda_i^{-1}d\lambda_i - d\lambda_k\lambda_k^{-1} - \lambda_k^{-1}d\lambda_i 
+ \lambda_k^{-1}d\lambda_k\lambda_k^{-1}\lambda_i\right)\tr(P_idP_kdP_k) 
\end{equation} 
where it is clear we lose nothing by restricting the sum to $i \neq k$. 

We now turn our attention to the second term appearing in~\eqref{eq: pullback of nu}: 
$$ 
\frac{\lambda_j\lambda_l\lambda_n}{\lambda_i\lambda_k\lambda_m}\tr(P_idP_j 
P_kdP_lP_mdP_n)
$$
Using the fact that $dP_i P_j + P_i dP_j = \delta_{ij}dP_i$ we 
can write this expression as 
\begin{align*} 
& \tr(\sum \frac{\lambda_j\lambda_l\lambda_n}{\lambda_i\lambda_k\lambda_m}
P_idP_jP_kdP_lP_mdP_n) \\ 
= & \tr(\sum \frac{\lambda_j\lambda_l\lambda_n}{\lambda_i\lambda_k\lambda_m}
(\delta_{mn}-P_n)P_i(\delta_{jk} - P_j)dP_kdP_ldP_m) \\ 
= & \tr(\sum \frac{\lambda_j\lambda_l\lambda_n}{\lambda_i\lambda_k\lambda_m}( 
P_i\delta_{mn}\delta_{jk} - \delta_{mn}\delta_{ij}P_j - \delta_{ni}\delta_{jk}P_i 
+ \delta_{ij}\delta_{in}P_j)dP_kdP_ldP_m) \\ 
= & \tr(\sum \left(\frac{\lambda_k}{\lambda_i} - \frac{\lambda_k}{\lambda_j} - 
\frac{\lambda_k}{\lambda_l} + \frac{\lambda_i\lambda_k}{\lambda_j\lambda_l}\right)
P_idP_jdP_kdP_l) 
\end{align*} 
We can finally write this as 
\begin{equation}
\label{eq: expression from Michael's office}  
\tr(\sum\left(\frac{\lambda_k(\lambda_l - \lambda_i)(\lambda_j - \lambda_i)}{ 
\lambda_i\lambda_j\lambda_l}\right)P_idP_jdP_kdP_l) 
\end{equation} 
This is a sum over four indices $i$, $j$, $k$ and $l$.  We make a series of observations to show that only certain 
combinations of these indices will make non-zero contributions to the sum.  This will allow us to eventually greatly 
simplify the sum.  
\begin{enumerate} 
\item We must have $i\neq j$ and $i\neq l$.  This is clear since $(\lambda_l - \lambda_i)(\lambda_j - \lambda_i)$ 
is a factor.  

\item The indices $i$, $j$, $k$ and $l$ cannot all be distinct.  If they were, then using the fact that $P_iP_j = 0$ and 
hence $dP_iP_j = - P_idP_j$ we see that the expression $\tr(P_idP_jdP_kdP_lP_i) = 0$.  Therefore some of $i$, $j$, $k$ and $l$ 
must be equal.  

\item No three of the indices $i$, $j$, $k$ and $l$ can be equal.  There are three possibilities to consider here: $i = j = k$, 
$i = j = l$ and $j = k = l$.  From the first observation above we can exclude the possibility that $i = j = k$ or $i = j = l$.  If 
$j = k = l$ then $\tr(P_idP_jdP_jdP_j) = - \tr(dP_idP_jdP_jP_jP_i) = 0$ after using the identities $P_idP_j = - dP_iP_j$ 
and $P_jdP_j + dP_jP_j = dP_j$.  

\item We conclude that two and only two of the indices $i$, $j$, $k$ and $l$ can be equal.  There are now six possibilities 
to consider (a) $i = j$, (b) $i = k$, (c) $i = l$, (d) $k = l$, (e) $j = l$, (f) $j = k$.  Of these six possibilities we can eliminate (a) 
and (b) immediately from the first observation.  The case where $i = k$ also gives no contribution since $\tr(P_dP_jdP_idP_l) = 
- \tr(dP_iP_jP_idP_ldP_i) = 0$.  The case where $j = l$ also gives no contribution for a similar sort of reason: $\tr(P_idP_jdP_kdP_j) = 
- \tr(dP_idP_jP_kP_jdP_i) = 0$.  
\end{enumerate}  
We conclude from this series of observations that we can restrict the sum to the two cases  $k = l$ or $j = k$ with all 
indices otherwise distinct.  By reindexing the sum we see that we can write~\eqref{eq: expression from Michael's office} as  
\begin{multline*}
\sum_{i, j, k\, \text{distinct}}\left(\frac{\lambda_j}{\lambda_i} - 1 - \frac{\lambda_j}{\lambda_k} + 
\frac{\lambda_i}{\lambda_k}\right)\tr(P_idP_jdP_jdP_k) \\ 
+ \sum_{i,j,k\, \text{distinct}}\left(\frac{\lambda_k}{\lambda_i} - \frac{\lambda_k}{\lambda_j} - 1 + 
\frac{\lambda_i}{\lambda_j}\right)\tr(P_idP_jdP_kdP_k) 
\end{multline*} 
Since $\tr(P_idP_jdP_kdP_k) = - \tr(P_jdP_kdP_kdP_i)$ for $i\neq j$ we can reindex the second sum to get 
$$ 
- \sum_{i,j,k\, \text{distinct}}\left(\frac{\lambda_j}{\lambda_k} - 
\frac{\lambda_j}{\lambda_i} - 1 + \frac{\lambda_k}{\lambda_i}\right)\tr(P_idP_jdP_jdP_k). 
$$ 
Adding these two sums gives 
$$ 
\sum_{i,j,k\, \text{distinct}}\left(2\frac{\lambda_j}{\lambda_i} - 2 
\frac{\lambda_j}{\lambda_k} + \frac{\lambda_i}{\lambda_k} - \frac{\lambda_k}{\lambda_i}\right) 
\tr(P_idP_jdP_jdP_k) 
$$ 
The three-forms $\tr(P_idP_jdP_jdP_k)$ for $\lambda_i$, $\lambda_j$ and 
$\lambda_k$ all distinct possess a cyclic symmetry
$$ 
\tr(P_idP_jdP_jdP_k) = \tr(P_kdP_idP_idP_j) = \tr(P_jdP_kdP_kdP_i).  
$$ 
Using this cyclic symmetry and re-indexing we get  
$$ 
3\sum_{i,j,k\, \text{distinct}}\left(\frac{\lambda_i}{\lambda_k} - 
\frac{\lambda_k}{\lambda_i}\right)\tr(P_kdP_idP_idP_j).  
$$ 
Clearly we can include the terms in this sum with $i=k$ with no effect, 
so that we can rewrite this as 
$$ 
3\sum_{i\neq j, j\neq k}\left(\frac{\lambda_i}{\lambda_k} - 
\frac{\lambda_k}{\lambda_i}\right)\tr(P_kdP_idP_idP_j).  
$$ 
Further, since $\tr(P_idP_jdP_jdP_j) = 0$ for $i\neq j$ we can write the 
sum as 
$$ 
3\sum_{j\neq k}\left(\frac{\lambda_i}{\lambda_k} - 
\frac{\lambda_k}{\lambda_i}\right)\tr(P_kdP_idP_idP_j).  
$$ 
This in turn can be written as 
$$ 
3\sum \left(\frac{\lambda_i}{\lambda_k} - 
\frac{\lambda_k}{\lambda_i}\right)\tr(P_kdP_idP_idP_j) - 3 \sum_{i\neq k}\left(\frac{\lambda_i}{\lambda_k} - 
\frac{\lambda_k}{\lambda_i}\right)\tr(P_kdP_idP_idP_k) 
$$ 
Again, using the fact that $\sum dP_j = 0$ we end up with 
$$ 
-3 \sum_{i\neq k} \left(\frac{\lambda_i}{\lambda_k} - 
\frac{\lambda_k}{\lambda_i}\right)\tr(P_kdP_idP_idP_k).  
$$ 
Reindexing we can write this as 
\begin{multline*}
-3\sum_{i\neq k}\left(\frac{\lambda_i}{\lambda_k} - 
\frac{\lambda_k}{\lambda_i}\right)\tr(P_kdP_idP_idP_j) + 
3\sum_{i\neq k}\frac{\lambda_i}{\lambda_k}\tr(P_idP_kdP_kdP_i) \\ 
= - 3\sum_{i\neq k} \frac{\lambda_i}{\lambda_k}\left( \tr(P_kdP_idP_idP_k) - \tr(P_idP_kdP_kdP_i)\right) 
\end{multline*} 
If $i\neq k$ then it is easy to show that 
$\tr(P_kdP_idP_idP_k) - \tr(P_idP_kdP_kdP_i) = \tr(dP_idP_kdP_i)$ and so we finally 
end up with, restoring the factor $-i/12\pi$, 
$$ 
-\frac{i}{4\pi}\sum_{i\neq k} \lambda_i\lambda_k^{-1}\tr(dP_idP_kdP_k),  
$$ 
since $\tr(dP_idP_kdP_k) = - \tr(dP_kdP_idP_i)$ if $i\neq k$.  
Therefore, combining this expression with~\eqref{eq: first part of p*nu} we end up with the 
following expression for $p^*(2\pi i\nu)$: 
\begin{multline*} 
\frac{i}{4\pi}\sum_{i \neq k}\left( 
\lambda_i^{-1}d\lambda_i - d\lambda_k\lambda_k^{-1} - \lambda_k^{-1}d\lambda_i 
+ \lambda_k^{-1}d\lambda_k\lambda_k^{-1}\lambda_i\right)\tr(P_idP_kdP_k) \\ 
- \frac{i}{4\pi}\sum_{i\neq k} \lambda_i\lambda_k^{-1}\tr(dP_idP_kdP_k).
\end{multline*} 
If we further pull this back to $\Omega^3(G/T\times Y_{T,\reg})$ then we find that it 
is equal to the expression~\eqref{eq: expression for p*df} we obtained for 
$ (p_Y)^*df$.  This completes the proof of Theorem~\ref{thm: main thm}.


\end{document}